\documentclass[]{interact}

\usepackage{multirow}
\usepackage{epstopdf}% To incorporate .eps illustrations using PDFLaTeX, etc.
\usepackage[caption=false]{subfig}% Support for small, `sub' figures and tables
\usepackage{xcolor}

\usepackage[numbers,sort&compress]{natbib}% Citation support using natbib.sty
\bibpunct[, ]{[}{]}{,}{n}{,}{,}% Citation support using natbib.sty
\renewcommand\bibfont{\fontsize{10}{12}\selectfont}% Bibliography support using natbib.sty

\theoremstyle{plain}% Theorem-like structures provided by amsthm.sty
\newtheorem{theorem}{Theorem}[section]
\newtheorem{lemma}[theorem]{Lemma}
\newtheorem{corollary}[theorem]{Corollary}

\theoremstyle{definition}

\theoremstyle{remark}
\newtheorem{remark}{Remark}

\def\ARMA{\textsc{arma}}
\def\AR{\textsc{ar}}

\def\ARCH{\textsc{arch}}
\def\GARCH{\textsc{garch}}

\newcommand{\tr}{{\rm ~tr}}
\newcommand{\VEC}{{\rm \,vec\,}}

\begin{document}

\articletype{ORIGINAL RESEARCH ARTICLE}% Specify the article type or omit as appropriate

\title{Bootstrapping a Powerful Mixed Portmanteau Test for Time Series}

\author{
\name{Esam Mahdi\textsuperscript{a}\thanks{CONTACT E Mahdi. Email: emahdi@qu.edu.qa} and Thomas J.\ Fisher\textsuperscript{b}}
\affil{\textsuperscript{a}Department of Statistical Sciences, University of Toronto, ON, Canada \textsuperscript{b}Department of Statistics, Miami University, Oxford, OH, USA}
}

\maketitle

\begin{abstract}
A new portmanteau test statistic is proposed for detecting nonlinearity in time series data. 
The new portmanteau statistic is calculated from the log of the determinant of a matrix comprised of the autocorrelations and cross-correlations of the residuals and squared residuals of a fitted time series.
The asymptotic distribution of the proposed test statistic is derived as a linear combination of chi-square distributed random variables and can be approximated by a gamma distribution. 
A bootstrapping approach is shown to be robust when distributional assumptions are relaxed.
The efficacy of the statistic is studied against linear and nonlinear dependency structures of some stationary time series models. 
It is shown that the new test can provide higher power than other tests in many situations. 
We demonstrate the advantages of the proposed test by investigating linear and nonlinear effects in {\color{blue} an economic series and} two environmental time series.
\end{abstract}

\begin{keywords}
Autocorrelation; Cross-correlation; Nonlinear time series; Portmanteau test.
\end{keywords}

%\linenumbers
\section{Introduction}

Whether in environmental or economic applications, in the modern practice of time series the detection of nonlinear dynamics, and the modeling thereof, is of fundamental importance.
Specifically, after a practitioner accounts for non-stationarity in a time series they will typically model any autocorrelation (or linear dynamics) as this is known to improve standard errors and forecasts. 
A logical step in this model fitting process is determining the adequacy of the fitted linear model.
In many cases, this determination is performed with a so-called Portmanteau statistic; see \cite{BoxPierce1970}.
Due to the skewness and heavy-tailed data we often incur in modern data analysis, computational methods are ften needed to relax distribution assumptions of many popular statistical approaches including the techniques for time series analysis.

Bootstrapping is a well-known computational method to approximate the variance, and other properties, of a sample statistic.
Since the seminal work of \cite{Efron:1979}, many variants of the bootstrap have been proposed and it is known to have many desirable properties \citep{Freedman:1981,EfronTibshirani:1994}.
In the field of time series it is important to retain any temporal structure in the bootstrap samples.
The block bootstrap \citep{Hall:1985,Kunsch:1989}, and stationary bootstrap \citep{PolitisRomano:1994}, were developed to work for time series.
Other variations exist \citep[][for example]{HallHorowitzJing:1995} and a review of the topic can be found in \cite{HardleHorowitzKreiss:2003}.
For regression applications with a heteroskedasticity, the Wild Bootstrap \citep{Wu:1986} was developed and its application has been well studied  \citep{Mammen:1993,MacKinnon:2006,MacKinnonEtAl:2021}.
This technique has been used in time series when {\color{blue} heteroskedasticity} is present \citep{Kim:2006,GospodinovTao:2011,ZhuEtAl:2020} and in other applications \citep[see][for example]{LeeBaek:2020}.
Recently, \cite{Lee:2016} proposed using the Wild bootstrap on the well-known Ljung-Box portmanteau test \cite{LjungBox1978} in time series and demonstrated it retains adequate type I error rates when heteroskedasticity is present.
A Randomly Weighted Bootstrap (RWB) (a variant of the Wild bootstrap proposed in \cite{JinYingWei:2001}) was proposed for several time series portmanteau statistics in the literature in the presence of heteroskedasticity in \cite{Zhu:2016}.
Recently, the RWB was used in the estimation of the Generalized Autoregressive Conditional Heterskedasticity (\GARCH) model in \cite{ZhuEtAl:2020}.

This article proposes a new portmanteau test that combines several existent results from the literature.
Utilizing the method of \cite{Zhu:2016} and \cite{ZhuEtAl:2020}, we use the RWB technique to approximate the distribution of the statistic under some fairly general scenarios (situations where the underlying stochastic process is heavy-tailed and skewed).
The article is organized as follows:\ Section \ref{sec:portmanteau.Statistics} provides a brief review of standard time series models and some popular portmanteau test statistics that have been used for detecting linear and nonlinear dependency in time series.
In Section \ref{sec:proposed.test} we propose a new portmanteau test statistic, derive its asymptotic distribution as a linear combination of chi-square random variables and discuss some of its properties. 
An extension of the RWB procedure is then outlined to approximate the distribution of the proposed statistic.
Section \ref{sec:computational.study} reports a Monte Carlo study comparing the empirical findings with the theoretical results and demonstrates that the empirical significance level of the proposed test statistic is accurately estimated by the percentiles of its asymptotic distribution.
Simulations also show that the power of the test is often higher than that of other test statistics. 
{\color{blue} Two illustrative applications are given in Section \ref{sec:application} to demonstrate the usefulness of the proposed test for real world datasets}. 
We end the article in Section \ref{sec:discussion} with some discussions on the advantages and limitations of the new statistic. 

%%%%%%%%%%%%%%%%%%%%%%%%%%%%%%%%%
%%%%%%%%%%%%%%%%%%%%%%%%%%%%%%%%%
%% Lit Review - Portmanteau
%%%%%%%%%%%%%%%%%%%%%%%%%%%%%%%%%
%%%%%%%%%%%%%%%%%%%%%%%%%%%%%%%%%

%%%%%=============================================================================%%%%%
\section{Time Series Modeling}\label{sec:portmanteau.Statistics}
%%%%%=============================================================================%%%%%

The autoregressive-moving average (\ARMA) model is arguably the most fundamental of all times series models.
An \ARMA$(p,q)$ for $n$ observations $z_1, z_2,\ldots ,z_n$ of a stationary mean $\mu$ time series can be expressed as
\begin{equation}\label{eq:arma}
	\Phi_p(B)(z_t-\mu) = \Theta_q(B)\varepsilon_t
\end{equation}
with
\begin{displaymath}
	\Phi_p(B) = 1 - \phi_1B - \phi_2B^2 - \ldots - \phi_pB^p,
\end{displaymath}
\begin{displaymath}
	\Theta_q(B) = 1 + \theta_1B + \theta_2B^2 + \ldots + \theta_qB^q,
\end{displaymath}
where $B$ is the backshift operator and the polynomials $\Phi_p(B)$ and $\Theta_q(B)$ are assumed to have all roots outside the unit circle on the complex plain and have no common roots.
The \emph{noise} sequence $\{\varepsilon_t\}$ may have further structure or be independent and identically distributed (iid) with mean 0 and constant variance, $\sigma^2>0$.

Let ${\boldsymbol\beta} = (\phi_1,\cdots,\phi_p,\theta_1,\cdots,\theta_q, \mu)$ denote the true parameter values and let $\hat{\boldsymbol\beta} = (\hat{\phi}_1,\cdots,\hat{\phi}_p,\hat{\theta}_1,\cdots,\hat{\theta}_q, \hat{\mu})$ denote the $\sqrt{n}$ consistent estimated values, so that the residuals $\hat{\varepsilon}_t$ denote the estimated values of $\varepsilon_t$ for $t = 1,\ldots, n$.
%A time series is considered \emph{linear} when it can be representaed by the General Linear Process, or when $\{z_t\}$ admits the {\it Moving-Average}, \MA, model  
%\begin{equation}\label{eq:linear}
%z_t-\mu_{t}= \Psi_\infty(B)\varepsilon_{t},
%\end{equation}
%where 
%$$\Psi_\infty(B)=\frac{\Theta_q(B)}{\Phi_p(B)}=1+\sum_{j=1}^{\infty}\psi_j(\bm{\zeta})B^j,~~\Phi_p(B)\ne 0~\hbox{for all}~|B|\leq 1,$$ 
%$\{\psi_j(\bm{\zeta})\}$ is a polynomial of parameters $\bm{\zeta}$ with absolutely summable sequence of weights, $\psi_j$. 
%When $\bm{\zeta}=(\phi_1,\cdots,\phi_p,\theta_1,\cdots,\theta_q)$, the \ARMA{} model in (\ref{eq:arma}) can be seen as a special case of %the linear representation given in (\ref{eq:linear}).
If the model in (\ref{eq:arma}) is correctly identified and the noise terms $\{\varepsilon_t\}$ are uncorrelated, then for all non-zero lag time $k$, the residual autocorrelation function, $\textrm{corr}\left(\varepsilon_t, \varepsilon_{t+k}\right)$, and the squared residual (or absolute-residual) autocorrelation function, $\textrm{corr}\left(\varepsilon_t^2, \varepsilon_{t+k}^2\right)$ (or $\textrm{corr}\left(\lvert\varepsilon_t\rvert, \lvert\varepsilon_{t+k}\rvert\right)$), should show no specific pattern and the correlation coefficient values should be approximately equal to zero. 
In addition, the cross-correlation function between the residuals and their squares, $\textrm{corr}\left(\varepsilon_t, \varepsilon_{t+k}^2\right)$, should be approximately uncorrelated with zero values.
On the other hand, if the model is not adequately identified, the autocorrelation may take on non-zero values.
Further, if there are nonlinear effects in the time series or if the residuals are not independent, these features may appear in the autocorrelation function of the squared (or the absolute) residuals or the cross-correlation of the residuals and their squares. 
The case for absolute residuals is beyond the scope of this article and we focus our attention on methods using the squared residuals.

Many nonlinear models have been proposed and can be used for analyzing nonlinear time series \cite[see Ch.\ 10 in][]{PTT2001}. 
For example, when the model is linear in mean but nonlinear in variance, \cite{Engle1982} proposed the Autoregressive Conditional heteroskedasticity, \ARCH, that is widely used for analyzing financial time series. 
This model was generalized by
\cite{Bollerslev1986}, the so-called \GARCH{} process. 
The innovations $\{\varepsilon_t\}$ in (\ref{eq:arma}) follows a \GARCH$(b, a)$ process if
\begin{equation}\label{eq:garch}
	\varepsilon_t = \xi_t\sigma_{t}, ~~~~ \sigma^2_{t} = \omega + \sum_{i=1}^b\alpha_i\varepsilon_{t-i}^2 + \sum_{j=1}^a \beta_j\sigma^2_{t-j},
\end{equation}
where $\xi_t$ is a sequence of iid random variables with a mean value of 0 and variance of 1, $\omega>0$, $\alpha_i\geq 0$ for $i=1,\ldots,a$ and $\alpha_i=0$ for all $i>a$, $\beta_j\geq 0$ for all $j=1,\ldots,b$ and $\beta_j=0$ for all $j>b$, and $\sum_{i=1}^{\max(b,a)} (\alpha_i+\beta_i) <1$. 
The summation on $\alpha_i\varepsilon_{t-i}^2$ terms comprises the original \ARCH{} model in \cite{Engle1982} with the second summation as the generalization.
Many variants of the \GARCH{} process have been proposed to model the dynamic behavior of conditional {\color{blue} heteroskedasticity} in real time series. 
Nice reviews of these models are available in \cite{Tsay2005} and \cite{Carmona2014}.

Typically, a practitioner may use a so-called portmanteau test to check the adequacy of a fitted model of form (\ref{eq:arma}) or for the presence of \GARCH-type effects in (\ref{eq:garch}).
Many authors have proposed such tests, including the commonly employed test statistics by \cite{BoxPierce1970,LjungBox1978}.
In, \cite{McLeodLi1983}, a portmanteau test for detecting nonlinearity (\textit{e.g.}, presence of \ARCH{} or \GARCH{}-effects) based on the squared residual autocorrelations is proposed. 
Several authors have improved on the portmanteau statistics by considering different functions of the autocorrelations of a fitted \ARMA{} model or autocorrelations of the squared residuals (see \cite{LiMak1994,PR2002,PR2006,PR2005,FisherGallagher2012}).
Simulation studies show that these statistics respond well in the detection of \ARCH{} models but tend to lack power compared to other types of nonlinear models that do not have the \ARCH{} effects. 

The idea of using empirical generalized correlations (correlation between $\hat{\varepsilon}_{t}^i$ and $\hat{\varepsilon}_{t+k}^j$, where $k$ is the lag time and $i,j$ are positive integers) to inspect for nonlinear dependency in time series models without considering the effects of the parameter estimation was introduced in \cite{LL1985,LL1987}.
Under the \ARMA{} assumptions in (\ref{eq:arma}),  
\cite{ZM2019} recently proposed portmanteau tests based on the generalized correlations and showed that the test based on the cross-correlation between the residuals and their squares can be more powerful than that of \cite{McLeodLi1983}, which is based on the autocorrelation of the squared residuals.

%%%%%%%%%%%%%%%%%%%%%%%%%%%%%%%%%%%%%%%%%%%%
\subsection{Portmanteau Review}
%%%%%%%%%%%%%%%%%%%%%%%%%%%%%%%%%%%%%%%%%%%%

%In this section, we review some well known portmanteau test statistics that have been used to detect the linearity and nonlinearity dependency in time series models.
Define $\hat{\rho}_{ij}(k)$ to be the correlation coefficient at lag time $k$ between $\hat{\varepsilon}_{t}^i$ and $\hat{\varepsilon}_{t+k}^j$, where we focus our attention on $i,j=1,2$. 
At lag time $k$, $\hat{\rho}_{11}(k)$ denotes the autocorrelation coefficient of the residuals,  $\hat{\rho}_{22}(k)$ denotes the autocorrelation coefficient of the squared residuals, and $\hat{\rho}_{12}(k)$ (or $\hat{\rho}_{21}(k)$) are the cross-correlation between the residuals and their squares at positive (or negative) lag $k$. 
Thus, $\hat{\rho}_{ij}(k)$ is given by  
\begin{multline}\label{eq:CC.sample}
	\hat{\rho}_{ij}(k)=\frac{\hat{\gamma}_{ij}(k)}{\sqrt{\hat{\gamma}_{ii}(0)}\sqrt{\hat{\gamma}_{jj}(0)}},\\~~\mbox{where}~~\hat{\gamma}_{ij}(k)=\frac{1}{n}\sum_{t=1}^{n-k}f_i(\hat{\varepsilon}_t)f_j(\hat{\varepsilon}_{t+k})~~\mbox{for}~~k=0,\pm 1, \pm 2,\cdots.
\end{multline}
Note that $\gamma_{ij}(k)=\gamma_{ji}(-k)$ for $k>0$ is the autocovariance (cross-covariance) at lag $k$ between the residuals to the power $i$ and the residuals to the power $j$ for $i,j=1,2$, for $f_i(x_{t})=x_{t}^i-n^{-1}\sum_{t=1}^{n}x_{t}^i$, for $i=1,2$.

Under the assumption that the data has been generated from an \ARMA{} process, \cite{BoxPierce1970} proposed to the time series literature the nominal portmanteau test in order to check the adequacy of the fitted model.
\cite{LjungBox1978} improved their test by utilizing a multiplicative factor on each squared autocorrelation term.
The two respective tests are 
\begin{equation}\label{eq:Ljungtest}
	\tilde{Q}_{11}=n\sum_{k=1}^{m}\hat{\rho}_{11}^2(k) ~~~\mbox{and}~~~Q_{11}=n(n+2)\sum_{k=1}^{m}(n-k)^{-1}\hat{\rho}_{11}^2(k),
\end{equation}
where $0<m<n/2$ is the maximum lag considered for significant autocorrelation.
Both $\tilde{Q}_{11}$ and $Q_{11}$ share the same asymptotic $\chi^2_{m-p-q}$ distribution but $Q_{11}$ generally has more power.

If the assumptions in (\ref{eq:arma}) are satisfied, \cite{McLeodLi1983} proposed a portmanteau test for detecting the presence of the \ARCH-effects, based on the autocorrelations of the squared residuals:\
\begin{equation}\label{eq:McLeodLitest}
	Q_{22}=n(n+2)\sum_{k=1}^{m}(n-k)^{-1}\hat{\rho}_{22}^2(k).
\end{equation}
\cite{McLeodLi1983} showed that the limiting distribution of $Q_{22}$ can be approximated by a chi-square distribution with $m$ degrees of freedom which different than $Q_{11}$ in (\ref{eq:Ljungtest}) since the limiting distribution of $Q_{22}$ does not depend on the order of the fitted \ARMA{} model.

Simulation studies show that the portmanteau statistics based on the squared residuals autocorrelations, such as $Q_{22}$, respond well to \ARCH{} models but tend to lack power in the presence of other types of nonlinear models. 
One possible reason for the lack of power could be the fact that these statistics ignore the generalized cross-correlation between the residuals to different powers; \textit{i.e.}, $\hat{\rho}_{ij}(k)$ in (\ref{eq:CC.sample}) for $i\neq j$. 
In this respect, \cite{LL1985,LL1987} proposed the idea of testing for nonlinearity in time series models using the cross-correlation between the residuals and squared residuals.
\cite{ZM2019} develop portmanteau tests to detect nonlinearity from stationary linear models, based on the generalized correlation.
Their test statistics are given by
\begin{equation}\label{eq:ZachariasMariantest1}
	Q_{ij}=n(n+2)\sum_{k=1}^{m}(n-k)^{-1}\hat{\rho}_{ij}^2(k),
\end{equation}
where $i,j=1,2,~i\ne j$.
These test can be seen as modified versions of the \cite{BoxPierce1970} and \cite{McLeodLi1983} tests that utilize the cross-correlation between the residuals and their squares. 
\cite{ZM2019} approximated the cross-correlation tests by $\chi_m^2$ and suggested that the tests based on the cross-correlations tend to be more powerful in detecting many types of nonlinearity compared to other statistics based on squared residual autocorrelations.

Many other portmanteau statistics have been developed for time series modeling.
\cite{PR2002} proposed a test to check the adequacy of the fitted \ARMA{} model based on the $m^\mathrm{th}$ root of the determinant of a $m^\mathrm{th}$ sample residual autocorrelations matrix.
They extended this statistic to test for nonlinearity by replacing the sample residual autocorrelations by the squared residual autocorrelations in the $m^\mathrm{th}$ autocorrelation matrix.
In \cite{PR2006}, they considered the log of the determinant of a $m^\mathrm{th}$ sample residual autocorrelations matrix.
These statistics are shown to be functions of the partial autocorrelation function, similar to \cite{Monti1994}.
\cite{MahdiMcLeod2012} extend the result to multivariate time series and \cite{FisherGallagher2012} uses the same matrix to {\color{blue}derive} a Weighted Ljung-Box and Weighted McLeod-Li test that are asymptotically similar to that of \cite{PR2002} and \cite{PR2006}.

\section{Proposed Test Statistic}\label{sec:proposed.test}
%%%%%=============================================================================%%%%%

Motivated by the results in \cite{LL1985,LL1987}, \cite{ZM2019}, and \cite{PR2002,PR2006}, we propose a new test for determining the adequacy of a fitted \ARMA{} model.
For a stationary time series, consider the block matrix of autocorrelations and cross-correlations of residuals and squared residuals, 
\begin{equation}\label{eq:GCMat}
	\hat{\bm{R}}(m) = \left[
	\begin{array}{c|c}
		\hat{\bm{R}}_{11}(m) & \hat{\bm{R}}_{12}(m)  \\ 
		\hline
		\hat{\bm{R}}_{12}^\prime(m)  & \hat{\bm{R}}_{22}(m) 
	\end{array}
	\right]_{2(m+1)\times 2(m+1)}
\end{equation}
\begin{multline*}
	\textrm{where}~~
	\hat{R}_{ii}(m) = \left[
	\begin{array}{cccc}
		1 & \hat{\rho}_{ii}(1) &\dots&\hat{\rho}_{ii}(m) \\
		\hat{\rho}_{ii}(-1) & 1 &\dots&\hat{\rho}_{ii}(m-1) \\
		\vdots & \vdots &\vdots&\vdots\\
		\hat{\rho}_{ii}(-m) & \hat{\rho}_{ii}(1-m) &\dots&1 
	\end{array}\right], ~~\textrm{for}~~ i=1,2 \\ \\
	\textrm{and}~~
	\hat{R}_{12}(m) = \left[\begin{array}{cccc}0 &\hat{\rho}_{12}(1)&\dots&\hat{\rho}_{12}(m)\\
		\hat{\rho}_{12}(-1)& 0 &\dots&\hat{\rho}_{12}(m-1)\\
		\vdots&\dots&\vdots&\vdots\\
		\hat{\rho}_{12}(-m)&\hat{\rho}_{12}(1-m)&\dots& 0
	\end{array}\right].
\end{multline*}
Note that $\hat{\bm{R}}_{21}(m)=\hat{\bm{R}}_{12}^\prime(m)$ is the matrix of cross-correlations between residuals and their squares, $\hat{\bm{R}}_{11}(m)$ is the residual autocorrelation Toeplitz matrix as defined in \cite{PR2002}, and $\hat{\bm{R}}_{22}(m)$ is the matrix of the autocorrelation coefficients 
based on the squared residuals.

Block matrices of this form have many desirable properties; \textit{e.g.}, this matrix is positive definite and $0<\lvert\hat{\bm{R}}(m)\rvert\leq 1$, where $\lvert\cdot\rvert$ denotes the determinant of the matrix.
Variants of $\hat{\bm{R}}(m)$ have been used to build other test statistics \cite[see][for example]{PR2002,MahdiMcLeod2012,RobbinsFisher2015} and is the foundation for our proposed portmanteau test:
\begin{equation}\label{eq:newtest}
	C_m=-\frac{n}{m}\log\lvert\hat{\bm{R}}(m)\rvert.
\end{equation}

Under the null hypothesis that an \ARMA{} model is adequately modeling the linear effects and no nonliner effects are present, $H_0:\bm{R}(m)=\bm{I}_{2(m+1)}$, where $\bm{I}_{2(m+1)}$ is the identity matrix of dimension $2(m+1)$.
That is, the sample autocorrelations of the residuals/squared residuals and the sample cross-correlations between the residuals and their squares will not significantly differ from zero, and $C_m$ will take on value near zero.
%That is, $\hat{\rho}_{ij}(k)\approx 0,~~(i,j =1,2)$ for $k=1,2,\cdots,m$ and $\hat{\rho}_{ij}(0)=1$ for $i,j=1,2$.
Under the alternative of an inadequate \ARMA{} model or the presence of nonlinear effects, $\bm{R}(m)$ will deviate from an identity, $\lvert\hat{\bm{R}}(m)\rvert<1$, and will approach 0 \cite[see][]{PR2002,RobbinsFisher2015}, thus $C_m$ will increase and the fitted linear model will be rejected as inadequate.

%%%%%%%%%%%%%%%%%%%%%%%%%%%%%%%%%%%%%%%%%%%%
\subsection{Distribution of Proposed Statistic}\label{sec:asymptotic.distribution}
%%%%%%%%%%%%%%%%%%%%%%%%%%%%%%%%%%%%%%%%%%%%

We now derive the asymptotic distribution of the proposed statistic $C_m$, while also discussing some complexities regarding it.
Through a sequence of Lemmas and Theorems we demonstrate that under the assumption of normality, $C_m$ is asymptotically distributed as a linear combination of $\chi^2$ random variables that can be approximated by a Gamma distribution.
We also discuss the distribution in the more general case.
In the derivation below $\bm{0}$ represents an appropriately sized vector or matrix of zeroes.

\begin{lemma}\label{lem:decomposeStat}
	The quantity $-n\log\lvert\hat{\bm{R}}(m)\rvert$ can asymptotically be decomposed into three components
	\begin{multline}\label{eq:detGCMat4}
		-n\log\lvert\hat{\bm{R}}(m)\rvert =-n\log\lvert\hat{\bm{R}}_{11}(m)\rvert -n \log\lvert\hat{\bm{R}}_{22}(m)\rvert \\
		+ n\left(\VEC \hat{\bm{R}}_{12}(m)\right)^{\prime}\left(\hat{\bm{R}}_{22}^{-1}(m)\otimes \hat{\bm{R}}_{11}^{-1}(m)\right)\VEC\hat{\bm{R}}_{12}(m),
	\end{multline}
\end{lemma}	
\begin{proof}
	The determinant of the block matrix $\hat{\bm{R}}(m)$ defined in (\ref{eq:GCMat}) is 
	$$\lvert\hat{\bm{R}}(m)\rvert=\lvert\hat{\bm{R}}_{11}(m)\rvert\times \lvert\hat{\bm{R}}_{22}(m)-\hat{\bm{R}}_{12}^{\prime}(m)\hat{\bm{R}}_{11}^{-1}(m)\hat{\bm{R}}_{12}(m)\rvert.$$
	Take the natural logarithm of the determinant and note that $\log(\lvert\bm{A}\rvert)=\tr\left(\log(\bm{A})\right)$, where $\tr(\bm{A})$ denotes the trace of matrix $\bm{A}$, thus
	\begin{displaymath}
		\log\lvert\hat{\bm{R}}(m)\rvert= \log\lvert\hat{\bm{R}}_{11}(m)\rvert + \tr\left(\log\left(\hat{\bm{R}}_{22}(m)-\hat{\bm{R}}_{12}^\prime(m)\hat{\bm{R}}_{11}^{-1}(m)\hat{\bm{R}}_{12}(m)\right)\right).
	\end{displaymath}
	
	Use a Taylor expansion \citep{Bhatia1997} of $\log\left(\hat{\bm{R}}_{22}(m)-\hat{\bm{R}}_{12}^\prime(m)\hat{\bm{R}}_{11}^{-1}(m)\hat{\bm{R}}_{12}(m)\right)$ and note that $$-\tr\left(\hat{\bm{R}}_{22}(m)\right)<\tr\left(\hat{\bm{R}}_{12}^\prime(m)\hat{\bm{R}}_{11}^{-1}(m)\hat{\bm{R}}_{12}(m)\right)<\tr\left(\hat{\bm{R}}_{22}(m)\right),$$ 
	this leads to
	\begin{multline}\label{eq:detGCMat2}
		\log\lvert\hat{\bm{R}}(m)\rvert\approx\log\lvert\hat{\bm{R}}_{11}(m)\rvert + \log\lvert\hat{\bm{R}}_{22}(m)\rvert \\ -\tr\left(\hat{\bm{R}}_{12}^\prime(m)\hat{\bm{R}}_{11}^{-1}(m)\hat{\bm{R}}_{12}(m)\hat{\bm{R}}_{22}^{-1}(m)\right).
	\end{multline}
	%where $R_{m}$ is the remainder term that is related to the dimension of the block matrix in (\ref{eq:GCMat}). 
	%Using the big-O definition and note that $\log(x)<0$ for $0<x<1$, one can show that $R_{m}\approx-m/n$.
	Recall that 
	$\tr\left(\bm{Z}^{\prime}\bm{BZC}\right)=(\VEC \bm{Z})^{\prime}(\bm{C}^{\prime}\otimes \bm{B})(\VEC \bm{Z})$, where $(\VEC\bm{Z})=[Z_{.1}:Z_{.2}:\cdots :Z_{.m}]^{\prime}$ is a column vector length $m^2$ formed by stacking the columns, $Z_{.j},j=1,2,\cdots,m$, of the $m\times m$ matrix $\bm{Z}$, and $\otimes$ is the Kronecker product of matrices \citep{Neudecker1969}.
	Thus (\ref{eq:detGCMat2}) becomes
	\begin{multline}\label{eq:detGCMat3}
		\log\lvert\hat{\bm{R}}(m)\rvert =\log\lvert\hat{\bm{R}}_{11}(m)\rvert + \log\lvert\hat{\bm{R}}_{22}(m)\rvert \\
		-\left(\VEC \hat{\bm{R}}_{12}(m)\right)^\prime\left(\hat{\bm{R}}_{22}^{-1}(m)\otimes \hat{\bm{R}}_{11}^{-1}(m)\right)\left(\VEC \hat{\bm{R}}_{12}(m)\right).
	\end{multline}
	Multiply (\ref{eq:detGCMat3}) by $-n$, the proposed statistic can be decompose into three components plus a constant term
	\begin{multline*}
		-n\log\lvert\hat{\bm{R}}(m)\rvert = -n\log\lvert\hat{\bm{R}}_{11}(m)\rvert -n \log\lvert\hat{\bm{R}}_{22}(m)\rvert\\
		+ n\left(\VEC \hat{\bm{R}}_{12}(m)\right)^{\prime}\left(\hat{\bm{R}}_{22}^{-1}(m)\otimes \hat{\bm{R}}_{11}^{-1}(m)\right)\left(\VEC \hat{\bm{R}}_{12}(m)\right).
	\end{multline*}
\end{proof}

\begin{lemma}\label{lemma:thirdComponentTrace}
	The third component in (\ref{eq:detGCMat4}) can be expressed as
	\begin{multline*}
		n\left(\VEC \hat{\bm{R}}_{12}(m)\right)^{\prime}\left(\hat{\bm{R}}_{22}^{-1}(m)\otimes \hat{\bm{R}}_{11}^{-1}(m)\right)\left(\VEC \hat{\bm{R}}_{12}(m)\right) \\
		=n\tr\left(\hat{\bm{R}}_{12}^{\prime}(m)\hat{\bm{R}}_{12}(m)\right). 
	\end{multline*}
\end{lemma}
\begin{proof}
	By $\VEC{(\bm{ABC})}=(\bm{C}^{\prime}\otimes \bm{A})(\VEC \bm{B})$ and 
	the fact  
	$(\VEC \bm{Z})^{\prime}(\bm{A}^{\prime}\otimes \bm{B}){\color{blue}(\VEC \bm{Z})}=\tr(\bm{AZ}^{\prime}\bm{BZ})$,
	we note that $\hat{\bm{R}}_{12}(m)=\bm{D}_{11}\hat{\bm{C}}_{12}(m)\bm{D}_{22}$, where $\bm{D}_{11}$ and $\bm{D}_{22}$ are, respectively, diagonal matrices with $i^\mathrm{th}$ diagonal elements $(\hat{\gamma}_{{11}_{ii}}(0))^{-1/2}$ and $(\hat{\gamma}_{{22}_{ii}}(0))^{-1/2}$, and $\hat{\bm{C}}_{12}(m)$ is the cross-covariance matrix of the residuals and their squares with order $m$.
	Now write the third component in (\ref{eq:detGCMat4}) as
	%\begin{equation*}
	{\small
		\begin{align*}
			n&\left(\VEC \hat{\bm{R}}_{12}(m)\right)^{\prime}\left(\hat{\bm{R}}_{22}^{-1}(m)\otimes \hat{\bm{R}}_{11}^{-1}(m)\right)\left(\VEC \hat{\bm{R}}_{12}(m)\right) \\
			&= n\left(\VEC \hat{\bm{C}}_{12}(m)\right)^{\prime}(\bm{D}_{22}\otimes \bm{D}_{11})\left(\hat{\bm{R}}_{22}^{-1}(m)\otimes \hat{\bm{R}}_{11}^{-1}(m)\right)(\bm{D}_{22}\otimes \bm{D}_{11})\VEC \hat{\bm{C}}_{12}(m)\\
			&= n\left(\VEC \hat{\bm{C}}_{12}(m)\right)^{\prime}\Big[\bm{D}_{22}\hat{\bm{R}}_{22}^{-1}(m)\bm{D}_{22}\otimes \bm{D}_{11}\hat{\bm{R}}_{11}^{-1}(m)\bm{D}_{11}\Big]\VEC \hat{\bm{C}}_{12}(m)\\
			&= n\left(\VEC \hat{\bm{C}}_{12}(m)\right)^{\prime}\left(\hat{\bm{C}}_{22}^{-1}(m)\otimes \hat{\bm{C}}_{11}^{-1}(m)\right)\VEC \hat{\bm{C}}_{12}(m)\\
			&= n\tr\left(\hat{\bm{C}}_{22}^{-1}(m) \hat{\bm{C}}_{12}^{\prime}(m)\hat{\bm{C}}_{11}^{-1}(m) \hat{\bm{C}}_{12}(m)\right)\\
			&= n\tr\left(\hat{\bm{C}}_{22}^{-1/2}(m) \hat{\bm{C}}_{12}^{\prime}(m)\hat{\bm{C}}_{11}^{-1/2}(m) \hat{\bm{C}}_{11}^{-1/2}(m) \hat{\bm{C}}_{12}(m)\hat{\bm{C}}_{22}^{-1/2}(m)\right)\\
			&= n\tr\left(\hat{\bm{R}}_{12}^{\prime}(m)\hat{\bm{R}}_{12}(m)\right).
	\end{align*}}
	%\end{equation*}
\end{proof}

Substituting the results of Lemma \ref{lemma:thirdComponentTrace} into (\ref{eq:detGCMat4}) of Lemma \ref{lem:decomposeStat} and multiply the results by the normalizing term $1/m$, the proposed statistic is asymptotically given by
\begin{equation}\label{eq:detGCMat5}
	\begin{split}
		C_m = -\frac{n}{m}\log\lvert\hat{\bm{R}}(m)\rvert=&-\frac{n}{m}\log\lvert\hat{\bm{R}}_{11}(m)\rvert \\
		&-\frac{n}{m} \log\lvert\hat{\bm{R}}_{22}(m)\rvert\\
		&+  \frac{n}{m}\tr\left(\hat{\bm{R}}_{12}^{\prime}(m)\hat{\bm{R}}_{12}(m)\right).
	\end{split}
\end{equation}
The proposed statistic has an interesting interpretation as it can be seen as an omnibus version of three existent tests.
The first component is asymptotically equivalent to the statistic proposed in \cite{PR2006}, $\tilde{D}_{11}$, which can be used to test for linear autocorrelation in the residuals. 
The second is also asymptotically equivalent to the test in \cite{PR2006},  $\tilde{D}_{22}$, which can be used to test for the nonlinearity, or heteroskedasticity, models (uncorrelated but not independent).
The third is asymptotically equivalent to a weighted variant \cite[in the vein of][]{FisherGallagher2012} of the tests 
proposed in \cite{ZM2019}, $Q_{12}$ and $Q_{21}$, which can be used to detect whether the cross correlations between the residuals and their squared values deviate from zero.

\begin{lemma}\label{lem:penaRodPart}
	The first and the second components in (\ref{eq:detGCMat5}) are given by
	\begin{equation}\label{eq:logpartialapproximate}
		-\frac{n\log\lvert\hat{\bm{R}}_{ii}(m)\rvert}{m}=-n\sum_{k=1}^{m}\frac{m+1-k}{m}\log(1-\hat{\pi}_{i,k}^2),~~i=1,2
	\end{equation}
	where $\hat{\pi}_{1,k}$ is the $k$th partial autocorrelation of residuals and 
	$\hat{\pi}_{2,k}$ is the $k$th partial autocorrelation of the squared residuals.
\end{lemma}
\begin{proof}
	The proof follows those presented in \cite{PR2002} and \cite{PR2005}.
\end{proof}

\begin{lemma}\label{lem:PsaradakisPart}
	The third component in (\ref{eq:detGCMat5}) can be written as {\color{blue} the sum of squares} of cross-correlation between residuals and their squares.
	\begin{equation}
		\small
		\frac{n}{m}\tr\left(\hat{\bm{R}}_{12}^{\prime}(m)\hat{\bm{R}}_{12}(m)\right)=
		\frac{n}{m}\sum_{\substack{k=-m,\\k\neq 0}}^{m}(m+1-|k|)\hat{\rho}_{12}^2(k).
	\end{equation}
\end{lemma}
\begin{proof}
	\begin{equation*}
		\small
		\begin{split}
			\frac{n}{m}\tr\left(\hat{\bm{R}}_{12}^{\prime}(m)\hat{\bm{R}}_{12}(m)\right)&=\frac{n}{m}\left(\VEC \hat{\bm{R}}_{12}(m)\right)^{\prime} \VEC \hat{\bm{R}}_{12}(m)\\
			&=\frac{n}{m}\Big[ (m+1)(0) + m[\hat{\rho}_{12}^2(-1)+\hat{\rho}_{12}^2(1)]\\
			&~~+(m-1)[\hat{\rho}_{12}^2(-2)+\hat{\rho}_{12}^2(2)]+
			\cdots+ [\hat{\rho}_{12}^2(-m)+\hat{\rho}_{12}^2(m)]\Big]\\
			&=\frac{n}{m}\sum_{\substack{k=-m,\\k\neq 0}}^{m}(m+1-\lvert k \rvert)\hat{\rho}_{12}^2(k).
		\end{split}
	\end{equation*}
\end{proof}

From the previous Lemmas, we rewrite the proposed test statistic as a linear combination of tests based on the partial autocorrelations of the residuals, partial autocorrelations of squared residuals, cross-correlation between the residuals and their squares at positive and negative lags as follows 
\begin{equation}\label{eq:detGCMat7}
	\begin{split}
		-\frac{n}{m}\log\lvert\hat{\bm{R}}(m)\rvert&= -\frac{n}{m}\sum_{k=1}^{m}(m+1-k)\log(1-\hat{\pi}_{1,k}^2)\\ &~~-\frac{n}{m}\sum_{k=1}^{m}(m+1-k)\log(1-\hat{\pi}_{2,k}^2)\\
		&~~+\frac{n}{m}\sum_{\substack{k=-m,\\k\neq 0}}^{m}(m+1-\lvert k \rvert)\hat{\rho}_{12}^2(k).
	\end{split}
\end{equation}

The marginal distribution of the first two components in (\ref{eq:detGCMat7}) can be found in \cite{PR2002}. 
Based on an application of the delta-method \citep[see][]{Monti1994,PR2002}, the distribution of the partial autocorrelations is asymptotically equivalent to the the autocorrelation and the first two components are asymptotically equivalent to the results in \cite{FisherGallagher2012}. For the third component, one can apply the methodology of \cite{RobbinsFisher2015} on the theoretical results of \cite{ZM2019} to find a similar marginal distribution. 
The asymptotic distribution of $C_m$ depends on the joint distribution of all terms $\hat{\rho}_{ij}(k)$ for $i,j=1,2$ and $k=1,\ldots,m$.

To derive the distribution of $C_m$ first note the equation in (\ref{eq:detGCMat7}) is asymptotically equivalent to a quadratic form $n\hat{\bm{r}}_m^\prime \bm{W}\hat{\bm{r}}_m$ where $\hat{\bm{r}}_m$ is a $4m$ vector comprised of autocorrelation and cross-correlation terms of residuals and their squares, and $\bm{W}$ is an appropriate diagonal matrix with elements corresponding to the weights associated with the components in (\ref{eq:detGCMat7}), given by 
$(m,m-1,\cdots,1,m,m-1,\cdots,1,m,m-1,\cdots,m,m-1,\cdots,1)/m$.

The distribution of statistics of the form $n\hat{\bm{r}}_m^\prime \bm{W}\hat{\bm{r}}_m$ is well understood when $\hat{\bm{r}}_m$ is asymptotically normal, see \cite{Box1954} and \cite{Satterthwaite1941,Satterthwaite1946}.
Thus, determining the distribution of $C_m$ is equivalent to determining the joint distribution of a vector mixed with autocorrelations of residuals, autocorrelation of squared residuals and cross correlation of residuals and their squares.
This problem is essentially an extension of \cite{WongLing:2005} and was recently studied in \cite{mahdi2021}.
In a general setting, knowledge of the joint third and fourth cumulants is necessary \citep[see][]{Brockwell1991} but in the case of iid Normal innovations, $\sqrt{n}\hat{\bm{r}}_m$ will be asymptotically normally distributed as follows.

%The statistic $C_m$ is a function of $\bm{r}_m$, which is a mix of partial autocorrelation and cross-correlation terms. 
%Using the arguments in can be shown using similar arguments as \cite{McLeodLi1983} and \cite{ZM2019}.
%Then, we will apply the result in \cite{LiMak1994} for an \ARCH model and modify this distribution to the case of \ARMA models with \ARCH errors.
%After that, the distribution of the proposed statistics follows using arguments in \cite{PR2002,PR2006,FisherGallagher2012,RobbinsFisher2015} based on Slutsky's theorem \cite[see][pp. 240-245]{Casella2001}.

\begin{theorem}\label{thm:propAsymNormal}
	If the assumptions in (\ref{eq:arma}) hold and the underlying stochastic process is normally distributed, then for any fixed integer $m<n$ the asymptotic distribution of 
	%\begin{equation*}
	%\sqrt{n}\hat{\bm{r}}_m = \sqrt{n}\Big(\hat{\rho}_{11}(1), \ldots, \hat{\rho}_{11}(m), %~\hat{\rho}_{22}(1), \ldots, \hat{\rho}_{22}(m),
	%~\hat{\rho}_{12}(1), \ldots, \hat{\rho}_{12}(m), 
	%~\hat{\rho}_{12}(-1), \ldots, \hat{\rho}_{12}(-m) \Big)^\prime 
	%\end{equation*}
	\begin{displaymath}
		\begin{split}
			\sqrt{n}\hat{\bm{r}}_m = \sqrt{n}\Big(\hat{\rho}_{11}(1), \ldots, \hat{\rho}_{11}(m),& ~\hat{\rho}_{22}(1), \ldots, \hat{\rho}_{22}(m),\\
			&\hat{\rho}_{12}(1), \ldots, \hat{\rho}_{12}(m), \hat{\rho}_{12}(-1), \ldots, \hat{\rho}_{12}(-m) \Big)^\prime 
		\end{split}
	\end{displaymath}
	is $\mathcal{N}(\bm{0}, {\boldsymbol\Sigma})$ where $\boldsymbol\Sigma$ is an $(4m)\times(4m)$ covariance matrix of the form
	\begin{equation}\label{eq:asymCovariance}
		{\boldsymbol\Sigma} = \left[\begin{array}{ccccc}
			\bm{I}_m - \bm{Q} & \bm{0} & \bm{0} & \bm{0} \\
			\bm{0} & \bm{I}_m & \bm{0} & \bm{0} \\
			\bm{0} & \bm{0} & \bm{I}_m & \bm{0} \\
			\bm{0} & \bm{0} & \bm{0} & \bm{I}_m 
		\end{array}\right]
	\end{equation}
	where $\bm{0}$ is $m\times m$ zero matrix, $\bm{Q}=\bm{X}_{m}\bm{V}^{-1}\bm{X}_{m}^\prime$ is an idempotent matrix with rank $p+q$, $\bm{V}$ is the information matrix for the parameters
	${\boldsymbol\beta}$ and $\bm{X}_m$ is an $m\times(p+q)$ matrix, with coefficients $\phi_i^\prime$ and $\theta_i^\prime$ defined by 
	$1/\phi(B)=\sum_{i=0}^{\infty}\phi_i^\prime B^i$ and $1/\theta(B)=\sum_{i=0}^{\infty}\theta_i^\prime B^i$ as defined in \cite[pp. 296-304]{Brockwell1991}.
\end{theorem}
\begin{proof}
	
	The proof may be established by straightforward calculation following similar arguments to that in \cite{WongLing:2005} where the results of \cite{BoxPierce1970} and  \cite{McLeodLi1983} are combined, but now include the results of \cite{ZM2019}.

\end{proof}

% \begin{remark}
	% Theorem \ref{thm:propAsymNormal} is presented based on the sample autocorrelations (largely for convenience of notation) whereas our proposed statistic is a function of the partial autocorrelation functions (\ref{eq:detGCMat7}). 
	% Note that the proof of Theorem \ref{thm:propAsymNormal} holds for partial autocorrelations by following the derivations in \cite{Monti1994} and \cite{PR2002}.
	% An application of the multivariate delta-method, see \cite[Appendix]{PR2002}, provides the result for $\log(1-\hat{\pi}_{i,k}^2)$.
	% Thus, under the adequacy of the fitted \ARMA models, our vector $\bm{r}_m$ is asymptotically normally distributed with the same covariance as that in (\ref{eq:asymCovariance}).
	% \end{remark}

\begin{theorem}\label{thm:proposedStatDistro}
	If the assumptions in Theorem \ref{thm:propAsymNormal} hold, then the asymptotic distribution of the proposed statistic is
	\begin{displaymath}
		C_m = -\frac{n}{m}\log\lvert\hat{\bm{R}}(m)\rvert \xrightarrow{D} \sum_{i=1}^{4m} \lambda_i \chi^2_{1,i}
	\end{displaymath}
	where $\chi^2_{1,i}$ are independent $\chi^2$ random variables and $\lambda_i$ are the eigenvalues of ${\boldsymbol\Sigma}\bm{W}$, with $\boldsymbol\Sigma$ defined in (\ref{eq:asymCovariance}) and $\bm{W}$ is a diagonal matrix with elements $(m,m-1,\cdots,1,m,m-1,\cdots,1,m,m-1,\cdots,1,m,m-1,\cdots,1)/m$.
\end{theorem}
\begin{proof}
	This results follows from the asymptotic normality in Theorem \ref{thm:propAsymNormal} along with the results in \cite{Box1954} and \cite{PR2002}.
\end{proof}

% \begin{remark}\label{remark:LiMak}
	% From the result in \citet{LiMak1994} for \ARCH$(b)$ models, $\sqrt{n}(\hat{\rho}_{22}(1), \ldots, \hat{\rho}_{22}(m))^\prime$ is asymptotically distributed as $\mathcal{N}(\bm{0}, {\boldsymbol I}_{m-(b+1)})$ where ${\boldsymbol I}_{m-(b+1)}$ is an $(m-(b+1))\times(m-(b+1))$ identity matrix.
	% Hence, for the \ARMA models with \ARCH errors, the submatrices ${\boldsymbol I}_m-{\boldsymbol Q}_m, {\boldsymbol I}_m, 1, {\boldsymbol I}_m, {\boldsymbol I}_m$ in the covariance block matrix given by (\ref{eq:asymCovariance}) will be replaced by ${\boldsymbol I}_m-{\boldsymbol Q}_m, {\boldsymbol I}_{m-(b+1)}, 1, {\boldsymbol I}_m, {\boldsymbol I}_m, respectively$. 
	% In this case, only a submatrix of diagonal ${\boldsymbol W}$ with weights $(m,m-1,\cdots,1,m,m-1,\cdots,(b+1),5(m+1),m,m-1,\cdots,1,m,m-1,\cdots,1)/(m+1)$ is used.
	% \end{remark}

\begin{corollary}\label{cor:proposedStatApprox}
	Under the assumptions of Theorem \ref{thm:propAsymNormal}, the asymptotic distribution of $C_m$ can be approximated by gamma distribution, $\Gamma(\alpha,\beta)$, where
	
	$$\alpha=\frac{3m[2(m+1)-(p+q)]^2}{4(m+1)(2m+1)-6m(p+q)},$$
	and
	$$\beta=\frac{4(m+1)(2m+1)-6m(p+q)}{3m[2(m+1) - (p+q)},$$
	where the distribution has a mean of $\alpha\beta=2(m+1)-(p+q)$ and a variance of $\alpha\beta^2=(4(m+1)(2m+1)-6m(p+q))/(3m)$.
\end{corollary} 
\begin{proof}
	Note that upper percentiles of the cumulative distribution of the form $\sum\lambda_i\chi_{1,i}^2$ can be approximated as $a\chi_{c}^2$, where the parameters $a$ and $c$ can be selected so that the mean and variance equal to those of exact distribution of $C_m$ \citep[see][]{Satterthwaite1941,Satterthwaite1946, Box1954}.
	Through cumulant matching arguments similar to \cite{PR2002} and \cite{FisherGallagher2012}, for large $m$ one can show that $a\chi_c^2$ is equivalent to a 
	gamma distribution with shape and scale parameters
	$$\alpha=K_1^2/K_2,~\hbox{and}~\beta=K_2/K_1,$$ 
	where
	\begin{equation*}
		K_1=\sum\lambda_i=\tr({\boldsymbol\Sigma}\bm{W})=2(m+1)-(p+q),
	\end{equation*}
	and
	\begin{equation*}
		K_2=2\sum\lambda_i^2=2\tr({\boldsymbol\Sigma}\bm{W})^2 = \frac{4(m+1)(2m+1)}{3m}-2(p+q).
	\end{equation*}
	From here, the result follows.
\end{proof}

%%%%%%%%%%%%%%%%%%%%%%%%%%%%%%%%%%%%%%%%%%%%
\subsection{Bootstrapping Algorithm}\label{sec:bootstrap}
%%%%%%%%%%%%%%%%%%%%%%%%%%%%%%%%%%%%%%%%%%%%

The asymptotic results presented in section \ref{sec:asymptotic.distribution} rely on iid normality of the innovations. 
In the modern practice of time series, violation of this assumption is common (consider the countless empirical examples where the series exhibits skewness and large tails). 
In these scenarios, the distribution of $C_m$ is more complex. 
Specifically, $\bm{r}_m$ will asymptotically be normally distributed but calculation of the covariance matrix requires derivation of the joint third and fourth cumulants of terms in $\bm{r}_m$ \citep[see Ch.\ 7.2 in][]{Brockwell1991}. 
To alleviate this complexity we adapt the RWB method proposed in \cite{Zhu:2016}, and recently used in \cite{ZhuEtAl:2020}, for use with our statistic. Specifically,
\begin{enumerate}
	\item Estimate the model from (\ref{eq:arma}) using least squares. From the residuals, compute $\hat{\rho}_{ij}(k)$ for $i,j=1,2$ and $k=1,\ldots,m$ and store in a vector $\hat{\bm{r}}_m$. 
	\item Generate a sequence of iid random variables $\bm{w}^* = \{w_1^*, w_2^*, \ldots, w_n^*\}$ independent of the data from a common distribution such that $P(w_i^*\geq 0)=1$ with mean and variance both equal to 1 (we use the standard Exponential distribution), and estimate the model (\ref{eq:arma}) using weighted least squares with weights $\bm{w}^*$ and compute the residuals, $\varepsilon_t^*$. 
	\item Calculate $\delta = \bm{W}\times\{\sqrt{n}(\hat{\bm{r}}_m^* - \hat{\bm{r}}_m)\}$, where $\hat{\bm{r}}_m^*$ is a length $4m$ vector with terms
	%    \begin{displaymath}
		%      \hat{\rho}_{ij}(k)^* = \frac{\sum_{t=k+1}^n w_t^*\left((\varepsilon_t^*)^i-E[(\varepsilon_t^*)^i]\right)\left((\varepsilon_{t-k}^*)^j-E[(\varepsilon_t^*)^j]\right)}{\sum_{t=1}^n \left((\varepsilon_t^*)^i-E[(\varepsilon_t^*)^i]\right)\left((\varepsilon_{t}^*)^j-E[(\varepsilon_t^*)^j]\right)}
		%    \end{displaymath}
	\begin{displaymath}
		\hat{\rho}_{ij}(k)^* = \frac{\sum_{t=k+1}^n w_t^*\left((\varepsilon_t^*)^i-E[(\varepsilon_t^*)^i]\right)\left((\varepsilon_{t-k}^*)^j-E[(\varepsilon_t^*)^j]\right)}{\sqrt{\sum_{t=1}^n \left((\varepsilon_t^*)^i-E[(\varepsilon_t^*)^i]\right)^2}\sqrt{\sum_{t=1}^n\left((\varepsilon_{t}^*)^j-E[(\varepsilon_t^*)^j]\right)^2}}
	\end{displaymath}
	with $E[(\varepsilon_t^*)^1]=0$ and $E[(\varepsilon_t^*)^2]=1$, and the matrix $\bm{W}$ is defined above. 
	\item Repeat steps 2 and 3 a large number of times, $B$ (typically $B\geq 500$), to obtain $\{\delta_{(1)}, \ldots, \delta_{(B)}\}$, and compute its covariance matrix and its associated eigenvalues $\hat{\lambda}_i^*$, for $i=1,\ldots, 4m$.
	\item Generate $N$ iid random samples $\{z_1^{(j)},\ldots,z_{4m}^{(j)}\}_{j=1}^N$, where $N$ is a large number (say $N=1000$), from a multivariate normal distribution with covariance $\bm{I}_{4m}$ and compute the sequence $\left\{K^{(j)}\right\}_{j=1}^N$ by
	\begin{displaymath}
		K^{(j)} = \sum_{i=1}^{4m} \hat{\lambda}_i^* \left(z_i^{(j)}\right)^2
	\end{displaymath}
	\item The sequence $\left\{K^{(j)}\right\}_{j=1}^N$ constitutes a bootstrapped sampling distribution for our proposed statistic $C_m$. Using the sample quantiles of $\left\{K^{(j)}\right\}_{j=1}^N$ we can determine critical values or we can approximate a $p$-value for $C_m$ by calculating $(\#(K^{(j)}>C_m)+1)/(N+1)$.
\end{enumerate}

The above algorithm is a logical extension of that proposed in \cite{Zhu:2016}, and a special case of that in \cite{ZhuEtAl:2020}.
The key to the algorithm is steps 3 and 4 where the covariance matrix of the vector $\bm{r}_m$ is approximated.

Unlike other bootstrapping methods the RWB approach does not require the practitioner to select a block length or similar parameters.
\cite{Zhu:2016} shows that the algorithm does not appear to be sensitive to the distribution of the weights (we found the standard exponential works fairly well).
Two parameters must be specified in the algorithm and are largely dependent on the computational resources available. 
In our simulations we use $B=2,000$ and $N=10,000$.

\begin{remark}
	The above algorithm can be modified for the other statistics discussed in this article, including $Q_{11}$, $Q_{22}$, $Q_{12}$, and $Q_{21}$ by only working with specific auto/cross-correlation terms in step 3.
\end{remark}

%%%%%%%%%%%%%%%%%%%%%%%%%%%%%%%%%
%%%%%%%%%%%%%%%%%%%%%%%%%%%%%%%%%
%% Computational Study
%%%%%%%%%%%%%%%%%%%%%%%%%%%%%%%%%
%%%%%%%%%%%%%%%%%%%%%%%%%%%%%%%%%

%%%%%=============================================================================%%%%%
\section{Computational Study}\label{sec:computational.study}
%%%%%=============================================================================%%%%%

We conduct a simulation study to investigate the appropriateness of the asymptotic distribution of the proposed test for different sample sizes and to compare its performance to the methods from the literature.
We also study the effects of skewness and excess kurtosis on the proposed method and demonstrate the bootstrapping algorithm in section \ref{sec:bootstrap} provides satisfactory results in approximating the distribution. 
Portmanteau statistics are known to be sensitive to the maximum lag, $m$, considered \citep[see][for a discussion]{GallagherFisher:2015}. 
For brevity, we limit our study to maximum lags $m=5$ and $m=10$.

We focus our attention on testing for the adequacy of a fitted \ARMA{} models. 
That is, our simulations consider the case of an underfit \ARMA{} model as well as the detection of nonlinear effects (\textit{e.g.}, \GARCH-type structures, or others) in the residual series.
We compare the proposed statistic, $C_m$, to that of $Q_{11}$ \cite{LjungBox1978}, $Q_{22}$ \cite{McLeodLi1983}, and the two tests $Q_{12}$ and $Q_{21}$ in \cite{ZM2019}.
The primary goals of our simulations are:\ to show that an omnibus, or mixed, test comprised of autocorrelations of residuals, their squares, and cross-correlation of the residuals and their squares, can gain in detection power of nonlinear models; and to show that the asymmetric structure of $C_m$ (where lag 1 terms appear $m$ times, lag 2 terms appear $m-1$ times, and so on) can improve power and help \emph{stabilize} the performance of a statistic across multiple lags.

To help facilitate the goals of our simulation, we also include the statistic
\begin{equation}\label{eq:comboStat}
	Q_{**} = n(n+2)\sum_{k=1}^m (n-k)^{-1}\left(\hat{\rho}_{11}^2(k) + \hat{\rho}_{22}^2(k) + \hat{\rho}_{12}^2(k) + \hat{\rho}_{21}^2(k)\right).
\end{equation}
Note that the statistic $Q_{**}$ is a combination of the $Q_{11}$, $Q_{22}$, $Q_{12}$ and $Q_{21}$, and following Theorem \ref{thm:propAsymNormal} will be approximately $\chi^2$ distributed with $4m-(p+q)$ degrees of freedom under the null hypothesis of (\ref{eq:arma}) and normal innovations.
The bootstrapping algorithm in section \ref{sec:bootstrap} can be utilized in the cases of non-normality where the $\bm{W}$ matrix in step 3 is replaced by a diagonal matrix with the Ljung-Box correction terms; \textit{i.e.}, $(n+2)/(n-k)$ for $k=1,\ldots,m$.

All numerical studies were conducted using the R software \citep{R2020} in a parallel framework with the \texttt{rugarch} package \citep{rugarch2020} for data generation.
This allows us to generate data with different nonlinear structures and under some fairly general distribution assumptions.
Source code is available in the supplementary material.

% We focus our attention on testing the adequacy of the fitted \ARMA($p,q$)+\GARCH($b,a$) models.
% We compare the proposed statistic, $C_m$ defined in (\ref{eq:newtest}), with the cross-correlation test statistics $(\tilde{Q}_{12},\tilde{Q}_{21}),(Q_{12},Q_{21})$ defined in (\ref{eq:ZachariasMariantest1}), the tests based on the autocorrelation of the residuals $Q_{11},Q_{11}^{w},M_{11}^{w},\tilde{D}_{11}$ defined in (\ref{eq:Ljungtest}), (\ref{eq:FisherGallaghertest1}), (\ref{eq:FisherGallaghertest2}), (\ref{eq:PR2002test}), and the tests based on the autocorrelation of the square-residuals
% $Q_{22},Q_{22}^{w},M_{22}^{w},\tilde{D}_{22}, L_b,L_{b}^{w}$  defined in (\ref{eq:McLeodLitest}), (\ref{eq:FisherGallaghertest1}), (\ref{eq:FisherGallaghertest2}), (\ref{eq:PR2002test}), (\ref{eq:FisherGallaghertest3}), respectively.
% In addition, we study the efficacy of the proposed statistic in detecting linearity and nonlinearity in several time series models. 
% All simulations were conducting use the R software \citep{R2020} and the source code is available in the supplementary material.

%%%%%%%%%%%%%%%%%%%%%%%%%%%%%%%%%%%%%%%%%%%%
\subsection{Studies on Empirical Size}\label{sec:sig.level}
%%%%%%%%%%%%%%%%%%%%%%%%%%%%%%%%%%%%%%%%%%%%

First we evaluate the empirical type I error rates of the {\color{blue} proposed} statistic, $C_m$, along with the others we will consider here, $Q_{11}$, $Q_{22}$, $Q_{12}$, $Q_{21}$ and $Q_{**}$, by calculating the rejection rate of the tests out of $1,000$ replications under the null hypothesis.
In Table \ref{tab:arNormalSize} we report the empirical size when the correct model is fit to a series of different sample sizes, $n=250$, $500$ and $1000$, generated by Gaussian \AR(1) with parameter $\phi=0.8$ and \AR(2) with $\phi_1=0.8$ and $\phi_2=-0.3$ processes.
In Table \ref{tab:arNormalSize}, generally we see all test report type I error rates within the acceptable range (3.7\% to 6.3\% based on Wald constructed 95\% acceptance regions).
Only in a few cases do we see rejection rates outside the acceptable range.
\begin{table}
\tbl{Empirical sizes at nominal rate of 5\% of $C_m$, $Q_{**}$, $Q_{11}$, $Q_{22}$, $Q_{12}$, and $Q_{21}$ under a Gaussian \AR(1) model with $\phi=0.8$ and Gaussian \AR(2) with $\phi_1=0.8$, $\phi_2=-0.3$ at different sample sizes and maximum lags $m$.}
{\begin{tabular}{lllccccccccccccccc}
				\multirow{2}{*}{$m$} &\multirow{2}{*}{$n$} & & \multicolumn{6}{c}{\AR(1) with $\phi=0.8$} & & & \multicolumn{6}{c}{\AR(2) with $\phi_1=0.8$, $\phi_2=-0.3$} \\
				\noalign{\smallskip}
				\cline{4-9}\cline{12-17}
				\noalign{\smallskip}
				&&&$C_m$ & $Q_{**}$ & $Q_{11}$ & $Q_{22}$ & $Q_{12}$ & $Q_{21}$ & & & $C_m$ & $Q_{**}$ & $Q_{11}$ & $Q_{22}$ & $Q_{12}$ & $Q_{21}$ \\
				\hline\noalign{\smallskip}
				&  250 & & 4.8 & 5.1 & 4.2 & 4.1 & 5.3 & 4.2 & & & 5.5 & 4.6 & 4.5 & 4.4 & 5.4 & 4.4\\
				5 &  500 & & 4.6 & 5.5 & 5.1 & 4.7 & 3.8 & 4.7 & & & 6.1 & 5.5 & 3.6 & 4.2 & 3.6 & 4.5 \\
				& 1000 & & 6.6 & 5.5 & 5.4 & 6.1 & 4.0 & 4.5 & & & 5.8 & 4.6 & 5.0 & 5.8 & 4.1 & 4.5 \\
				\noalign{\smallskip}\hline
				\noalign{\smallskip}
				&  250 & & 5.6 & 5.3 & 4.3 & 5.3 & 4.7 & 4.6 & & & 5.5 & 6.7 & 4.7 & 4.8 & 4.7 & 4.6 \\
				10 &  500 &  &5.8 & 5.7 & 5.5 & 5.8 & 4.8 & 4.1 & & & 6.0 & 5.6 & 5.5 & 6.2 & 5.0 & 4.7 \\
				& 1000 & & 5.6 & 5.3 & 5.0 & 5.2 & 3.8 & 4.4 & & & 5.5 & 4.7 & 4.9 & 6.1 & 3.8 & 4.9 \\
				\hline
		\end{tabular}}
\label{tab:arNormalSize}
\end{table}

The presented theoretical findings are based on the assumption of an underlying Gaussian process as the distribution of $C_m$ (and $Q_{**}$) in a more general setting is more complicated.
Table \ref{tab:arSkewHeavySize} displays the empirical size at the nominal rate of 5\% when utilizing the asymptotic distribution under the same \AR(1) process where the innovation are not normally distributed.
On the left side we report the type I error rates when the innovations are generated from the Skewed Normal distribution such that the skewness is approximately 0.56 (this corresponds to skewness parameter 1.5 in the \texttt{rdist} function in the \texttt{rugarch} package).
On the right side the innovations are from the Students' $t$ distribution where the excess kurtosis is 1 (shape parameter 10 in the \texttt{rdist} function).

\begin{table}
\tbl{Empirical sizes at nominal rate of 5\% of $C_m$, $Q_{**}$, $Q_{11}$, $Q_{22}$, $Q_{12}$, and $Q_{21}$ under a \AR(1) process with $\phi=0.8$ when the innovations are generated from a Skewed Normal or Students' $t$ distribution, at different sample sizes and maximum lags $m$.}
{\begin{tabular}{lllccccccccccccccc}
				\multirow{2}{*}{$m$} &\multirow{2}{*}{$n$} & & \multicolumn{6}{c}{Skewed Normal Innovations} & & & \multicolumn{6}{c}{Students' $t$ Innovations} \\
				\noalign{\smallskip}
				\cline{4-9}\cline{12-17}
				\noalign{\smallskip}
				&&&$C_m$ & $Q_{**}$ & $Q_{11}$ & $Q_{22}$ & $Q_{12}$ & $Q_{21}$ & & & $C_m$ & $Q_{**}$ & $Q_{11}$ & $Q_{22}$ & $Q_{12}$ & $Q_{21}$ \\
				\hline\noalign{\smallskip}
				&  250 & & 6.9 & 6.5 & 5.6 & 4.7 & 2.8 & 3.9 & & & 6.6 & 6.2 & 5.4 & 5.9 & 4.6 & 4.9\\
				5 &  500 & & 6.9 & 6.6 & 4.7 & 3.8 & 5.0 & 5.6 & & & 6.7 & 6.1 & 4.4 & 5.7 & 5.6 & 4.8 \\
				& 1000 & & 6.0 & 5.7 & 6.1 & 4.9 & 3.8 & 4.1 & & & 5.7 & 5.0 & 3.9 & 5.3 & 5.5 & 3.9 \\
				\noalign{\smallskip}\hline
				\noalign{\smallskip}
				&  250 & & 7.6 & 7.0 & 4.1 & 6.1 & 4.1 & 4.8 & & & 6.8 & 5.7 & 5.2 & 5.8 & 5.1 & 4.5 \\
				10 &  500 & & 8.0 & 6.6 & 5.3 & 4.7 & 5.2 & 4.9 & & & 6.0 & 5.4 & 4.1 & 4.9 & 5.0 & 3.8 \\
				& 1000 & & 7.1 & 5.9 & 4.3 & 5.2 & 4.6 & 4.3 & & & 6.2 & 6.4 & 4.9 & 5.4 & 4.2 & 5.4 \\
				\hline
		\end{tabular}}
\label{tab:arSkewHeavySize}
\end{table}

Table \ref{tab:arSkewHeavySize} shows the proposed statistic $C_m$ and the combination statistic $Q_{**}$, which use a combination of autocorrelations of residuals, squared residuals and the cross-correlation of the residuals and their squares, begins to report inflated type I error rates in the presence of skewness.
This phenenomen appears more problematic at the larger maximum lag of $m=10$.
In the case of heavy tails, we only see a moderate increase in type I error rates.

In table \ref{tab:arSSTDSize} we consider the robustness of the statistics when the RWB algorithm is utilized.
Data is generated from the same \AR(1) process above with innovations from the Skewed Students' $t$ distribution such that the skewness is approximately 0.85 and the excess kurtosis is 1.73, thus the underlying innovations come from a distribution with both heavy-tails and skewness.
The table reports the empirical rejection rates at a nominal rate of 5\% when the asymptotic distribution is utilized and when the RWB algorithm in Section \ref{sec:bootstrap} is used.
For comparison, we also include the RWB-based type I error rates for the statistics $Q_{11}$, $Q_{22}$, $Q_{12}$ and $Q_{21}$.
\begin{table}
\tbl{Empirical sizes at nominal rate of 5\% of $C_m$, $Q_{**}$, $Q_{11}$, $Q_{22}$, $Q_{12}$, and $Q_{21}$ under an \AR(1) model with $\phi=0.8$ and at different sample sizes and maximum lags $m$ when the innovations are generaed from a Skewed Students' $t$ distribution. Results with the asymptotic distribution and the randomly weighted bootstrap algorithm are presented.}
{\begin{tabular}{lllccccccccccccccc}
				\multirow{2}{*}{$m$} &\multirow{2}{*}{$n$} & & \multicolumn{6}{c}{Based on asymptotic distribution} & & & \multicolumn{6}{c}{Based on RWB algorithm} \\
				\noalign{\smallskip}
				\cline{4-9}\cline{12-17}
				\noalign{\smallskip}
				&&&$C_m$ & $Q_{**}$ & $Q_{11}$ & $Q_{22}$ & $Q_{12}$ & $Q_{21}$ & & & $C_m$ & $Q_{**}$ & $Q_{11}$ & $Q_{22}$ & $Q_{12}$ & $Q_{21}$ \\
				\hline\noalign{\smallskip}
				&  250 & & 7.2 & 6.7 & 6.0 & 4.6 & 5.2 & 3.8 & & & 2.7 & 1.8 & 5.5 & 2.4 & 3.8 & 3.0\\
				5 &  500 & & 7.8 & 7.1 & 4.7 & 5.3 & 3.6 & 5.0 & & & 2.7 & 1.7 & 4.5 & 3.7 & 3.6 & 1.5  \\
				& 1000 & & 7.4 & 8.0 & 5.0 & 5.7 & 4.0 & 4.5 & & & 2.9 & 2.6 & 4.8 & 3.1 & 4.1 & 2.9 \\
				\noalign{\smallskip}\hline
				\noalign{\smallskip}
				&  250 & & 8.9 & 9.3 & 6.1 & 5.0 & 4.1 & 5.2 & & & 1.9 & 1.2 & 4.9 & 1.7 & 3.1 & 1.9 \\
				10 &  500 & & 8.0 & 7.0 & 7.0 & 6.0 & 4.1 & 4.0 & & & 1.6 & 1.3 & 6.4 & 1.4 & 3.1 & 1.4  \\
				& 1000 & & 8.5 & 8.3 & 4.7 & 5.7 & 5.3 & 3.9 & & & 1.4 & 1.9 & 4.2 & 2.4 & 4.2 & 1.6 \\
				\hline
		\end{tabular}}
\label{tab:arSSTDSize}
\end{table}

Table \ref{tab:arSSTDSize} shows that the proposed statistic and $Q_{**}$, both which are asymptotically equivalent to linear combinations of $Q_{11}$, $Q_{22}$, $Q_{12}$ and $Q_{21}$, report inflated type I error rates when utilizing the asymptotic distribution.
This is due to the joint third and fourth cumulants.
As before, we also note the increased type I error rate appears larger for the larger maximum lag of $m=10$.
However, when utilizing the RWB algorithm to approximate the distribution, the type I errors do not exceed the nominal level for the proposed statistic.
Overall we see most statistics report conservative type I error rates, which is generally preferred in practice compared to inflated type I errors.

%%%%%%%%%%%%%%%%%%%%%%%%%%%%%%%%%%%%%%%%%%%%
\subsection{\label{sec:power.study}Power Studies}
%%%%%%%%%%%%%%%%%%%%%%%%%%%%%%%%%%%%%%%%%%%%

We now consider the empirical power of the proposed method and compare it to some of the statistics in the literature.
The first two simulation scenarios are structured to demonstrate the proposed statistic $C_m$ provides comparable power to one of the statistics that should be more powerful.
Specifically, we generate data from a Gaussian \ARMA(1,1) model with $\phi_1=0.8$ and $\theta_1=0.3$ but only an \AR(1) model is fit; thus we've intentionally underfit the autocorrelation in the series.
We would expect temporal correlation in the residuals, thus the Ljung-Box test $Q_{11}$ should detect the underfit.
In the second scenario, the data follows a Gaussian \AR(1)+\ARCH(1) process with $\phi_1=0.8$, $\alpha_1=0.4$ and $\omega=1$ but only an \AR(1) is fit.
Here, we would expect substantial temporal correlation in the squares of the residuals and for the McLeod-Li test $Q_{22}$ to be quite powerful.
We highlight the most powerful statistic in \textbf{boldface}.

\begin{table}
\tbl{Empirical power at nominal rate of 1\% of $C_m$, $Q_{**}$, $Q_{11}$, $Q_{22}$, $Q_{12}$, and $Q_{21}$ under two alternatives, a Gaussian \ARMA(1,1) model with $\phi_1=0.8$ and $\theta_1=0.3$ and a Gaussian \AR(1)+\ARCH(1) with $\phi_1=0.8$, $\omega=1$ and $\alpha_1=0.4$ where the process is underfit with an \AR(1) at different sample sizes and maximum lags $m$.}
{\begin{tabular}{lllccccccccccccccc}
				\multirow{2}{*}{$m$} &\multirow{2}{*}{$n$} & & \multicolumn{6}{c}{Underfit \ARMA(1,1)} & & & \multicolumn{6}{c}{Underfit \AR(1)+\ARCH(1)} \\
				\noalign{\smallskip}
				\cline{4-9}\cline{12-17}
				\noalign{\smallskip}
				&&&$C_m$ & $Q_{**}$ & $Q_{11}$ & $Q_{22}$ & $Q_{12}$ & $Q_{21}$ & & & $C_m$ & $Q_{**}$ & $Q_{11}$ & $Q_{22}$ & $Q_{12}$ & $Q_{21}$ \\
				\hline\noalign{\smallskip}
				&  250 & & 63.0 & 34.0 & \textbf{73.1} & 2.4 & 2.1 & 1.5 & & &  \textbf{77.8} & 65.5 & 5.3 & 75.2 & 10.2 & 11.6\\
				5 &  500 & & 97.1 & 83.5 & \textbf{98.6} & 4.7 & 1.1 & 1.3 & & & \textbf{98.1} & 94.4 & 4.8 & 98.0 & 13.5 & 19.0 \\
				& 1000 & & \textbf{100.0} & 99.9 & \textbf{100.0} & 7.1 & 1.2 & 1.3 & & & \textbf{100.0} & 99.8 & 4.0 & \textbf{100.0} & 15.3 & 23.8  \\
				\noalign{\smallskip}\hline
				\noalign{\smallskip}
				&  250 & & 51.4 & 24.2 & \textbf{57.8} & 2.1 & 1.8 & 1.8 & & & \textbf{70.3} & 52.2 & 3.2 & 65.6 & 8.1 & 8.9 \\
				10 &  500 & & 93.9 & 67.2 & \textbf{96.2} & 2.9 & 1.5 & 1.2 & & & 96.3 & 88.1 & 3.8 & \textbf{96.6} & 11.3 & 16.0  \\
				& 1000 & & \textbf{100.0} & 99.4 & \textbf{100.0} & 5.4 & 1.5 & 1.4 & & & 99.8 & 99.5 & 3.3 & \textbf{99.9} & 11.6 & 18.7 \\
				\hline
		\end{tabular}}
\label{tab:armaArchPower}
\end{table}

Table \ref{tab:armaArchPower} reports the empirical power at the 1\% significance level.
For the underfit \ARMA(1,1) process we see that the Ljung-Box statistic is most powerful with the proposed statistic $C_m$ providing comparable power.
Given in this scenario all temporal structure should be present in the residuals, and since the proposed $C_m$ is comprised of terms involving the residuals and their squares, it is not overly surprising the traditional Ljung-Box test is more powerful.
In the second scenario we see either the proposed $C_m$ or McLeod-Li $Q_{22}$ is most powerful.
We also note that the cross-correlations of the residuals and their squares (the components of $Q_{12}$ and $Q_{21}$) provide some detection (although not particularly strong) of the underlying \ARCH{} process.
This partially explains why the proposed method is most powerful in a few scenarios -- it combines elements of $Q_{22}$ with $Q_{12}$ and $Q_{21}$.
Lastly, we note that the weighted statistic $C_m$ is more powerful than its non-weighted counterpart $Q_{**}$ in all scenarios and that all statistics demonstrate a reduction in power as the lag increases.
This phenomenon is well-known in the literature \cite[see][for example]{GallagherFisher:2015} but note that the decrease in power appears less for the proposed statistic, $C_m$.

The next study considers detecting some nonlinear processes studied in the literature. 
Data is generated from the following five processes:
%\begin{equation}\label{keenan.models}
%\begin{split}
\begin{flalign}\label{eq:nonlinearModels}
	\nonumber\hbox{M1},~~~~z_t&=\varepsilon_t-0.3 \varepsilon_{t-1} +0.2 \varepsilon_{t-2} +0.4 \varepsilon_{t}\varepsilon_{t-2}-0.25\varepsilon_{t-2}^2,&&\\
	\nonumber\hbox{M2},~~~~z_t&=0.4z_{t-1} -0.3z_{t-2}+0.5z_{t-1}\varepsilon_{t-1}+\varepsilon_{t},&&\\
	\nonumber\hbox{M3},~~~~z_t&=0.4z_{t-1} -0.3z_{t-2}+0.5z_{t-1}\varepsilon_{t-1}+0.8\varepsilon_{t-1}+\varepsilon_{t},&&\\
	\hbox{M4},~~~~z_t&=0.5-(0.4-0.4\varepsilon_{t-1})z_{t-1}+\varepsilon_{t},&&\\
	\nonumber\hbox{M5},~~~~z_t&=0.8\varepsilon_{t-2}^2+\varepsilon_t.
\end{flalign}
%\end{split}
%\end{equation}
where $\varepsilon_t$'s are a sequence of independent and identically distributed innovations.
The first three models are analyzed by \cite{Keenan1985} (see also \cite{PR2002,PR2006}), whereas the other models are studied in \cite{ZM2019}.
The nonlinear process is then fit with an \AR($p$) where $p$ is selected using the Akaike Information Criterion (AIC) and statistics are calculated based on the estimated residuals.
To study the robustness of the statistics we generate the innovations from a Skewed Students' $t$-distribution and utilize the RWB algorithm for all the statistics.

\begin{table}
\tbl{Empirical power at nominal rate of 5\% of $C_m$, $Q_{**}$, $Q_{11}$, $Q_{22}$, $Q_{12}$, and $Q_{21}$ under various nonlinear alternatives where the process is fit with an \AR($p$) with $p$ selected based on AIC, at different sample sizes and maximum lags $m$.}
{\begin{tabular}{llccccccccccccccc}
				\multirow{2}{*}{Model} & & \multicolumn{6}{c}{Lag $m=5$} & & & \multicolumn{6}{c}{Lag $m=10$} \\
				\noalign{\smallskip}
				\cline{3-8}\cline{11-16}
				\noalign{\smallskip}
				&&$C_m$ & $Q_{**}$ & $Q_{11}$ & $Q_{22}$ & $Q_{12}$ & $Q_{21}$ & & & $C_m$ & $Q_{**}$ & $Q_{11}$ & $Q_{22}$ & $Q_{12}$ & $Q_{21}$ \\
				\hline\noalign{\smallskip}
				& &\multicolumn{14}{c}{\underline{$n=250$}}\\ 
				\noalign{\smallskip}
				M1 & & 76.2 & 64.1 & 0.3 & 29.3 & 47.1 & \textbf{84.6} & & & 72.5 & 48.5 & 0.2 & 22.5 & 37.7 & \textbf{73.1} \\
				M2 & & \textbf{84.6} & 76.9 & 0.6 & 74.7 & 36.1 & 74.8 & & & \textbf{85.0} & 72.6 & 1.1 & 73.1 & 32.0 & 70.1  \\
				M3 & & 89.4 & 84.8 & 0.9 & 79.4 & 13.3 & \textbf{90.6} & & & \textbf{88.5} & 80.2 & 1.0 & 76.8 & 11.9 & 88.4  \\
				M4 & & \textbf{79.5} & 57.6 & 0.0 & 57.2 & 55.0 & 30.1 & & & \textbf{79.1} & 51.5 & 0.1 & 55.0 & 49.9 & 24.9  \\
				M5 & & \textbf{30.0} & 20.2 & 0.6 & 20.8 & 11.1 & 22.7 & & & \textbf{27.6} & 11.9 & 0.8 & 16.6 & 8.2 & 15.2 \\
				\hline\noalign{\smallskip}
				& &\multicolumn{14}{c}{\underline{$n=500$}}\\ 
				\noalign{\smallskip}
				M1 & & 88.3 & 82.9 & 0.0 & 54.4 & 75.4 & \textbf{93.2 }& & & 88.7 & 78.6 & 0.2 & 47.3 & 70.2 & \textbf{90.9} \\
				M2 & & \textbf{88.6} & 84.1 & 0.4 & 81.3 & 36.4 & 86.9 & & & \textbf{88.5} & 81.6 & 0.6 & 80.2 & 32.7 & 84.8 \\
				M3 & & 93.0 & 89.4 & 1.3 & 84.4 & 11.4 & \textbf{94.9} & & & 92.5 & 87.0 & 0.6 & 83.4 & 10.6 & \textbf{93.4}  \\
				M4 & & \textbf{84.6} & 70.9 & 0.0 & 73.1 & 56.8 & 59.5 & & & \textbf{85.0} & 69.1 & 0.2 & 72.5 & 54.6 & 52.5 \\
				M5 & & \textbf{50.3} & 40.7 & 0.6 & 33.7 & 8.9 & 49.0 & & & \textbf{43.2} & 27.4 & 0.4 & 27.4 & 6.0 & 35.6  \\
				\hline
		\end{tabular}}
\label{tab:nonlinearPower}
\end{table}

Table \ref{tab:nonlinearPower} reports the empirical power of the statistics using the RWB algorithm at the nominal level of 5\%.
The results demonstrate the proposed statistics is comparable, or more powerful, than those studied from the literature.
In particular, we draw attention to the results for models M2, M4 and M5.
There, each of $Q_{22}$, $Q_{12}$ and $Q_{21}$ provides substantial power but $C_m$, which uses information from all three of those statistics, provides a notable increase in power.
We also see the expected behavior of an increase in power as $n$ increases and a general decrease in power as the lag $m$ increases.

Our last simulation considers the potential increase in power offered by our statistic, and the potential gains of using an omnibus statistic. 
Consider a modification of M4 above:\
\begin{equation}\label{eq:model4}
	z_t =0.5-(0.4-\delta\varepsilon_{t-1})z_{t-1}+\varepsilon_{t}\\
\end{equation}
where $\delta$ can be considered a perturbation parameter that controls the amount of \textit{nonlinearity}.
When $\delta=0$ we have a simple \AR(1) process, as $\delta$ increases the process becomes increasingly nonlinear.

\begin{figure}[htp!]
	\begin{center}
		\includegraphics[width=0.75\textwidth]{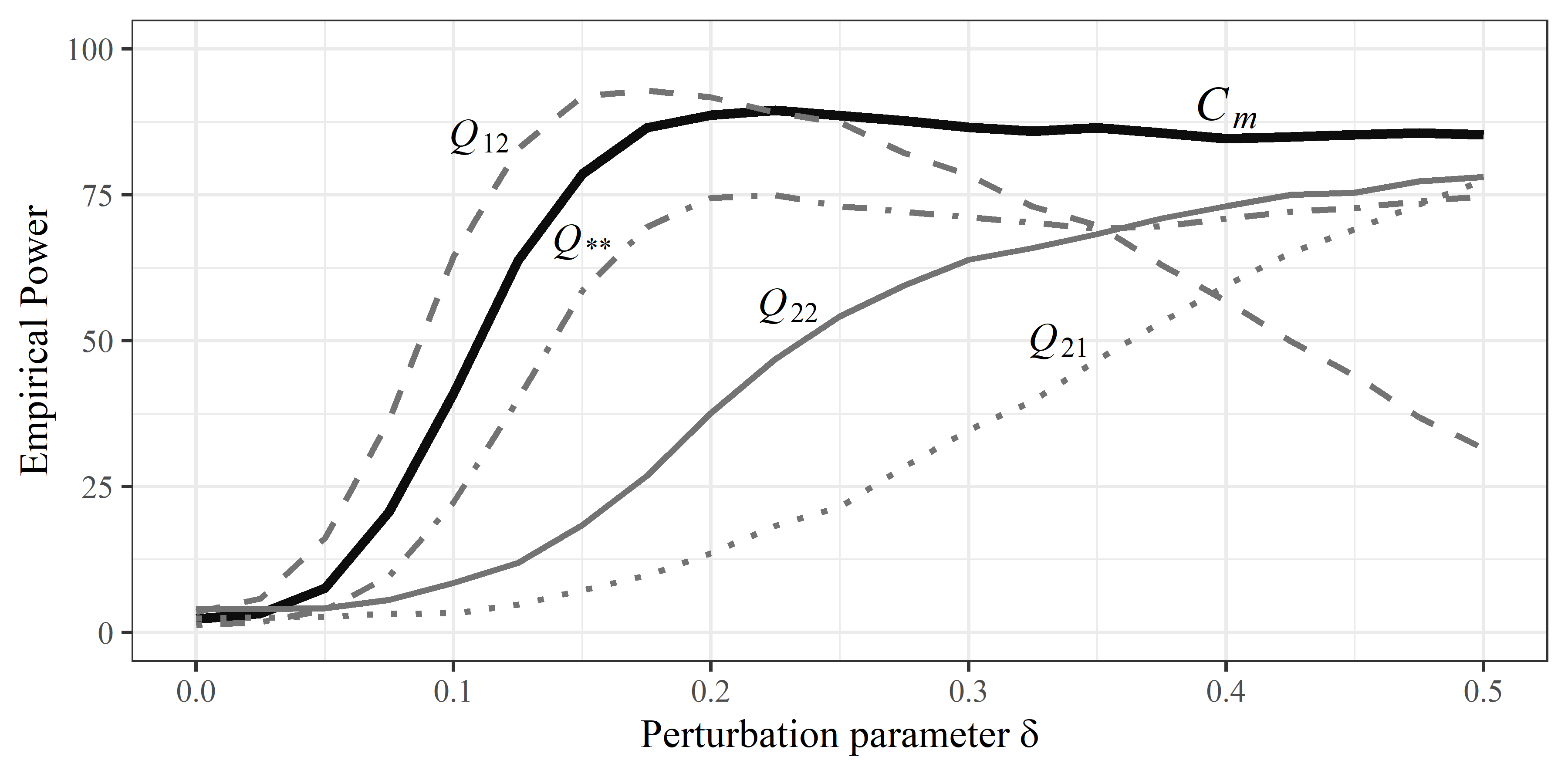}
	\end{center}
	\caption{{\color{blue}Powers of $C_m$ (solid black line), $Q_{**}$ (dotted dashed gray line), $Q_{22}$ (solid gray line), $Q_{12}$ (dashed gray line) and $Q_{21}$ (dotted gray line) at lag $m=5$ in detecting the nonlinear structure in (\ref{eq:model4}) for a 5\% test when data of size 500 are generated and an \AR($p$) is fit to the data with $p$ selected from AIC.}}
	\label{fig:nonlinearPower}
\end{figure}

Data is generated using model (\ref{eq:model4}) for $n=500$ with Skewed Students' $t$-distributed innovations and values of $\delta$ ranging from 0 to 0.5.
An \AR($p$) is fit to the data, where $p$ is selected from AIC (generally we may expect AIC to select $p=1$ to model the autoregressive part of the process). 
Figure \ref{fig:nonlinearPower} displays the empirical power of the statistics $C_m$, {\color{blue} $Q_{**}$,} $Q_{22}$ (second most powerful statistic in Table \ref{tab:nonlinearPower}), $Q_{12}$ and $Q_{21}$.
The distribution is approximated through the RWB algorithm and the empirical power is calculated based on 1,000 realizations at each $\delta$ value.

Figure \ref{fig:nonlinearPower} shows a very interesting power plot. 
At smaller $\delta$ values, we see that the statistic $Q_{12}$ is most powerful followed closely by the proposed statistics $C_m$. 
Yet, as $\delta$ increases the power of $Q_{12}$ starts to decrease while the statistics $Q_{22}$ and $Q_{21}$ gain power. 
The proposed statistic, $C_m$, which is asymptotically a convolution of the other statistics, effectively combines all the information and is most powerful starting around $\delta=0.225$.
{\color{blue} The power of $Q_{**}$ follows the pattern of $C_m$ but is generally less powerful, perhaps suggesting most of the remaining correlation is at lower lags \citep[see][]{GallagherFisher:2015}.}

\begin{figure}[htp!]
	\begin{center}
		\includegraphics[width=0.85\textwidth]{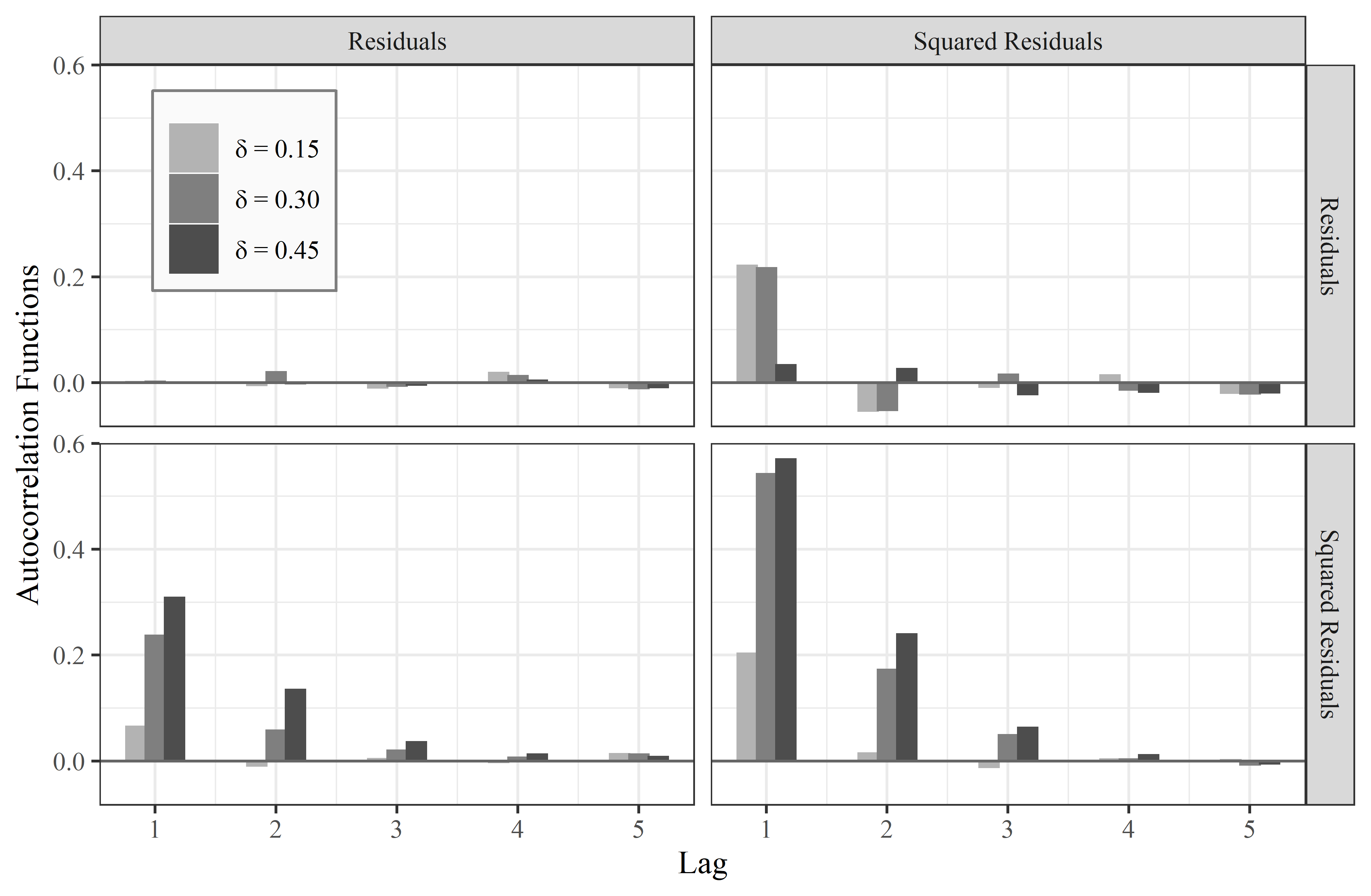}
	\end{center}
	\caption{Correlograms of the residuals and their squares from simulated datasets from nonlinear model (\ref{eq:model4}) at three different $\delta$ values when an \AR($p$) is fit to the data, $p$ is selected by AIC.}
	\label{fig:nonlinearPowerExploration}
\end{figure}

To understand the behavior of the statistics in Figure \ref{fig:nonlinearPower}, we conduct some follow up simulations to study the nature of the autocorrelation and cross-correlation of the residuals and their squares from the above study.
Three series of length $n=5,000$ were simulated from model (\ref{eq:model4}) (we chose a large $n$ to ensure some level of consistency in the estimation of the residual autocorrelation and cross-correlation functions), one each for $\delta=0.15$, $0.30$ and $0.45$.
Each of the three nonlinear series was fit with an \AR($p$) where $p$ is chosen from AIC (in our simulations $p=2$, $1$, $2$ were selected, respectively) and Figure \ref{fig:nonlinearPowerExploration} displays the sample autocorrelation of the residuals and their squares (on the diagonal), and the cross correlation of the residuals and their squares (on the off diagonal).

Figure \ref{fig:nonlinearPowerExploration} shows that for $\delta=0.15$ the strongest residual correlations occur at $\hat{\varepsilon}_{12}(1)$ and $\hat{\varepsilon}_{12}(2)$ which correspond to components of the $Q_{12}$ statistic. 
There is also meaningful correlation in the $\hat{\varepsilon}_{22}(1)$ term (corresponding to the moderate power in the $Q_{22}$ near $\delta=0.15$).
At $\delta=0.30$, there is meaningful correlation in $\hat{\varepsilon}_{12}(k)$, $\hat{\varepsilon}_{21}(k)$ , and $\hat{\varepsilon}_{22}(k)$ terms up to $k=3$, which helps demonstrate why all the statistics demonstrate some detection power, and that $C_m$ is most powerful.
At $\delta=0.45$ little correlation remains in the components of $Q_{12}$, which explains its decrease in power, while there is substantial correlation in the components of $Q_{22}$ (and $Q_{21}$ to a lesser extent).
This corresponds closely to the results of Table \ref{tab:nonlinearPower} which shows $C_m$ as most powerful, followed by $Q_{22}$ and then $Q_{21}$ slighlty more powerful than $Q_{12}$.

This example provides an interesting case study on the usefulness of an omnibus statistic, such as $C_m$.
Without oracle type knowledge, a practitioner would be unable to rely on a single statistic such at $Q_{12}$, $Q_{21}$ or $Q_{22}$, to detect nonlinearity in the residual series.
The proposed omnibus statistic $C_m$ can encapsulate all the relevant information, all while weighting the components in a way known to increase power as seen compared to the studied $Q_{**}$ statistic \citep[see][]{GallagherFisher:2015}.

In conclusion, the simulations demonstrate the proposed statistic $C_m$ can attain good power (at least comparable to other methods, if not better) compared to many of the proposed statistics in the literature.
The simulations demonstrate that by using all the information contained in the autocorrelations of the residuals, autocorrelation of the squared residuals and cross correlation of the residuals and their squares, one can attain more power in detecting nonlinear effects than any statistic based on just one measure, all while retaining adequate type I error rates and providing comparable power in detecting underfit linear effects.

%%%%%%%%%%%%%%%%%%%%%%%%%%%%%%%%%
%%%%%%%%%%%%%%%%%%%%%%%%%%%%%%%%%
%% Illustrative Applications
%%%%%%%%%%%%%%%%%%%%%%%%%%%%%%%%%
%%%%%%%%%%%%%%%%%%%%%%%%%%%%%%%%%

%%%%%=============================================================================%%%%
\section{\label{sec:application}Illustrative Applications}
%%%%%=============================================================================%%%%

We demonstrate the usefulness of the proposed test for detecting nonlinear processes in {\color{blue} an economic series and} some environmental data recently studied in the literature.

{\color{blue}
\subsection{\label{sec:econApplication} Crude Oil Prices}

Consider a short study on the daily West Texas Intermediate (WTI) Crude Oil Prices \citep{WTIdata}, in U.S. dollars per barrel, from September 01, 2019 through July 20, 2022 obtained using the \texttt{tidyquant} package \citep{tidyquant2022}. 
This time frame encompasses all market days beginning roughly six months before the onset of the \emph{lock downs} due to the SARS-CoV-2 pandemic and the economic turbulence that has occurred since, and results in a length $n=723$ series. 
The changes in daily crude prices ($\nabla\textrm{WTI}$) are seen in Figure \ref{fig:wtiChanges} along with the normal QQ-plot of the daily returns showing the data is heavy-tailed.
\begin{figure}[htp!]
	\begin{center}
		\includegraphics[width=0.85\textwidth]{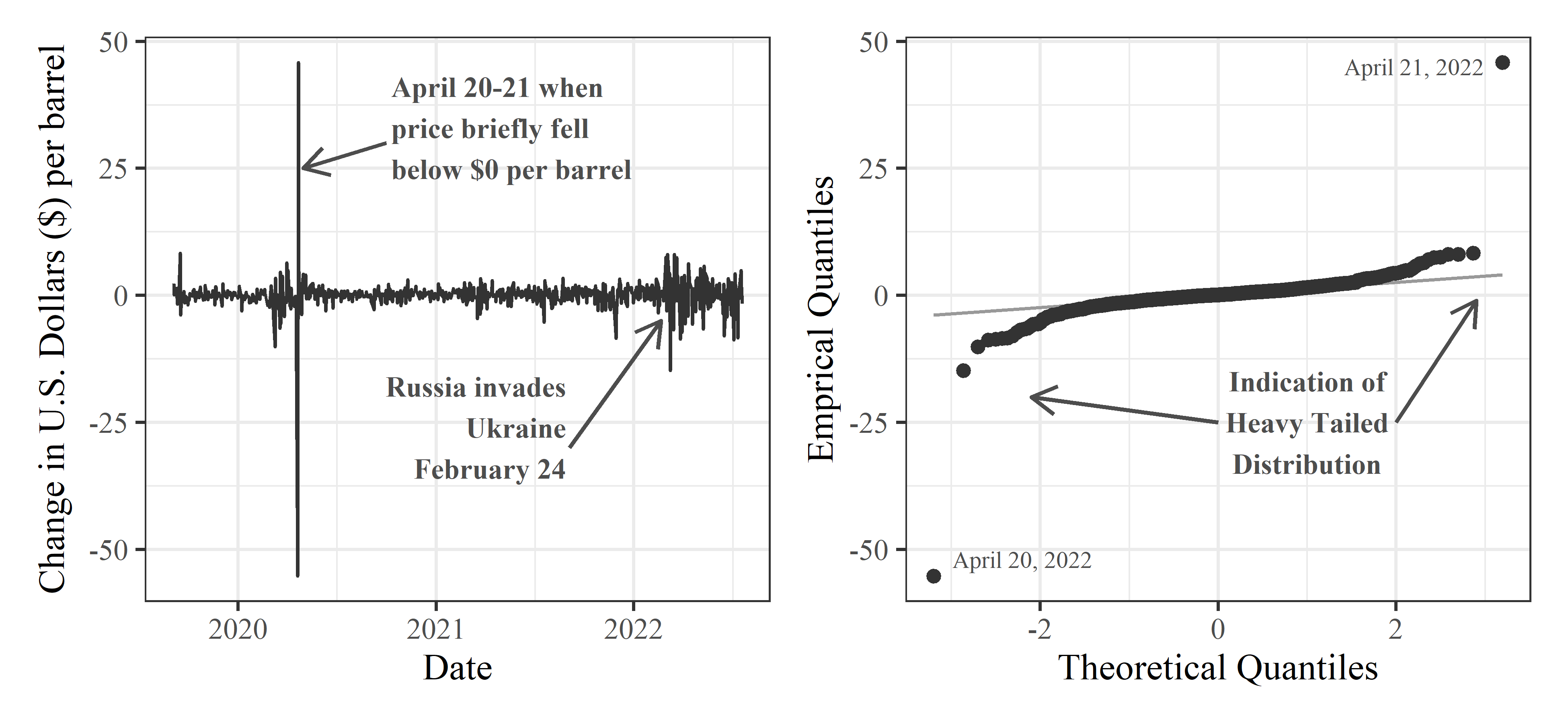}
	\end{center}
	\caption{{\color{blue}Daily returns for the West Texas Intermediate Crude Oil Prices (USD per barrel) from September 2019 through July 2022 (left) and normal QQ-Plot (right) demonstrating a heavy tailed distribution.}}
	\label{fig:wtiChanges}
\end{figure}

Figure \ref{fig:wtiChanges} shows changes in daily crude prices ($\nabla\textrm{WTI}$) are reasonably stationary but that the distribution of changes exhibit heavy tails.
The series exhibit a decaying autocorrelation function (not shown) and AIC suggest an \AR(3) will model the linear dependency in the data.
The fitted \AR(3) model has parameters
\begin{displaymath}
	\phi_1=-0.3144, ~\phi_2=-0.1443, ~\phi_3=-0.0733, ~~\sigma^2=10.44 
\end{displaymath}
and the autocorrelation functions of the resulting residual series and squared residuals series can be seen in Figure \ref{fig:wtiResACFs} (note the off-diagonal terms are the cross correlations of the residuals and their squares).
\begin{figure}[htp!]
	\begin{center}
		\includegraphics[width=0.85\textwidth]{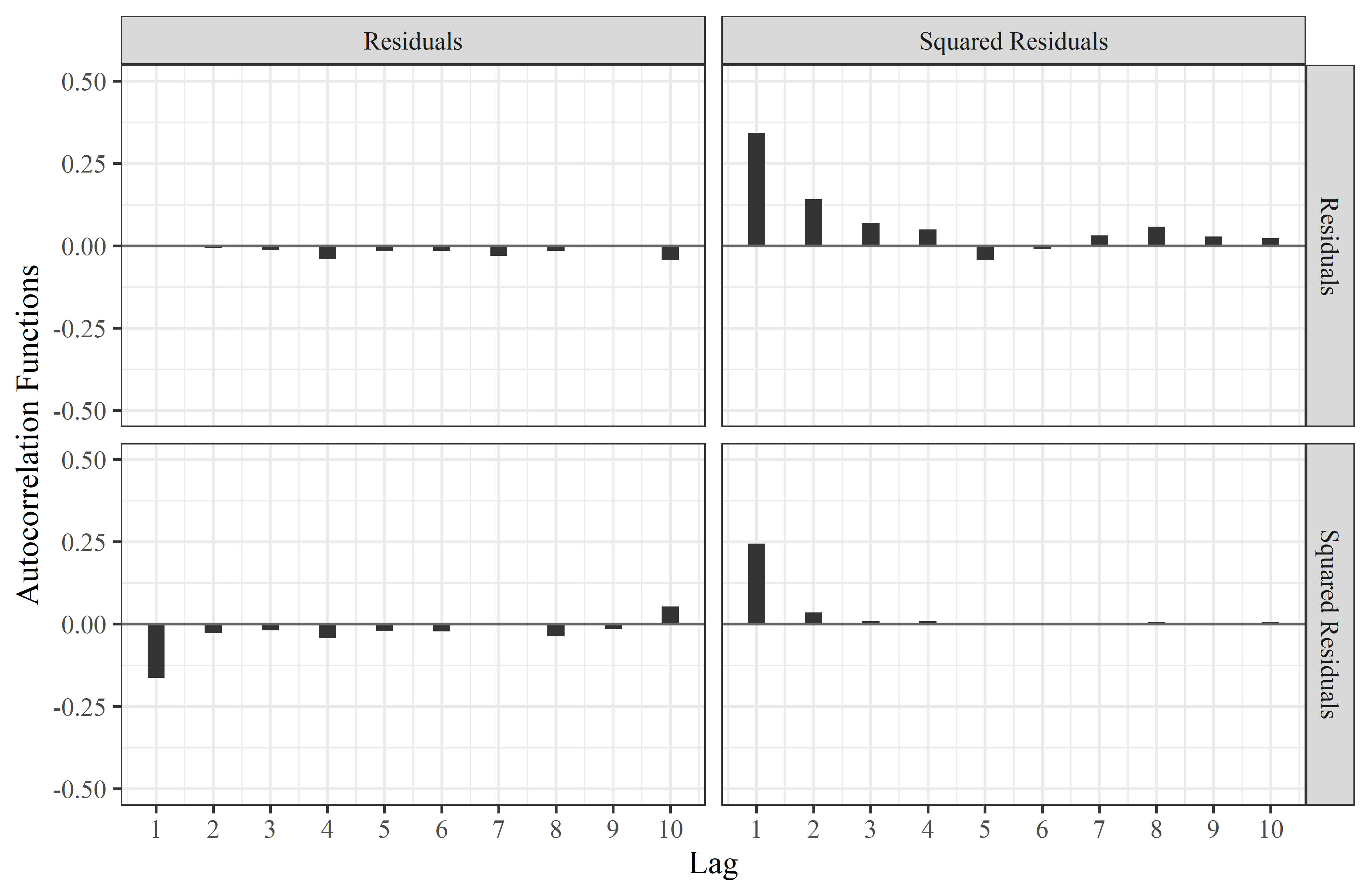}
	\end{center}
	\caption{{\color{blue}Correlograms of the residuals and their squares from the \AR(3) model fit to the daily returns for the WTI Cude Oil Price series.}}
	\label{fig:wtiResACFs}
\end{figure}
There we see no meaningful temporal correlation in the residuals (indicating an \AR(3) adequately models the linear dependence structure). 
There appears to be meaningful correlation in the $\hat{\rho}_{12}(k)$ terms and measurable correlation in $\hat{\rho}_{21}(1)$ and $\hat{\rho}_{22}(1)$.
Given the non-normality of the data we apply the RWB algorithm to compute the various test statistics at lags $m=5$, $m=10$ and $m=20$ with the results provided in Table \ref{tab:wtiPvals} (we use $B=10,000$ and $N=10,000$ in this application).
\begin{table}
\tbl{{\color{blue}The statistics and associated $p$-values, based on the RWB algorithm, of the portmanteau tests when an \AR(3) is fit to the daily returns for the WTI series.}}
{\begin{tabular}{llrrccrrccrr}
					&& \multicolumn{2}{c}{Lag $m=5$} & & & \multicolumn{2}{c}{Lag $m=10$} & & & \multicolumn{2}{c}{Lag $m=20$}\\
					\noalign{\smallskip}
					\cline{3-4}\cline{7-8}\cline{11-12}
					\noalign{\smallskip}
					& &Stat.\ & RWB & & & Stat.\ & RWB & & & Stat.\ & RWB \\
					& &       & $p$-value & & & &  $p$-value & & & &  $p$-value \\
					\hline\noalign{\smallskip}
					\noalign{\smallskip}
					$C_m$ & & 170.02 & 0.008  & & & 176.26 & 0.009 & & & 196.19 & 0.010 \\
					$Q_{**}$ & & 174.22 & 0.009 & & & 184.41 & 0.012 & & & 226.27 & 0.019\\
					$Q_{11}$ & & 1.56 & 0.617 & & & 3.80 & 0.886  & & & 25.74 & 0.309 \\
					$Q_{22}$ & & 44.64 & 0.002 & & & 44.70 & 0.002 & & & 44.92 & 0.002 \\
					$Q_{12}$ & & 106.26 & 0.022 & & & 110.54 & 0.023 & & & 116.14 & 0.022 \\
					$Q_{21}$ & & 21.77 & 0.005 & & & 25.36 & 0.029  & & & 39.48 & 0.119 \\
					\hline
			\end{tabular}}
\label{tab:wtiPvals}
\end{table}
Not surprisingly given Figure \ref{fig:wtiResACFs}, we see the Ljung-Box test $Q_{11}$ confirms we have adequately modeled the linear process. 
The tests of Pasaradakis-V\'{a}vra \cite{ZM2019}, $Q_{12}$ and $Q{21}$, both reject at the smaller lags but provide differing results at the higher lag of $m=20$, and the McLeod-Li \cite{McLeodLi1983} test, $Q_{22}$ also detects nonlinearity.
Both omnibus test, $C_m$ and $Q_{**}$, reject the null hypothesis that a linear model is adequate.

We note that the $p$-value of $Q_{**}$ more than doubles when the lag increases from $m=5$ to $m=20$ while that of the $C_m$ is relatively constant.
In general, a practitioner must choose the lag $m$ and many of the portmanteau test of the Ljung-Box form (including $Q_{**}$) are known to be sensitive to the lag \citep[see][]{ GallagherFisher:2015}.
To study these effects we compute the $p$-value of $C_m$ and $Q_{**}$ using the RWB algorithm at lags $m=4,\ldots,40$ and display them in Figure \ref{fig:wtiPvaluePlot}.
\begin{figure}[htp!]
	\begin{center}
		\includegraphics[width=0.55\textwidth]{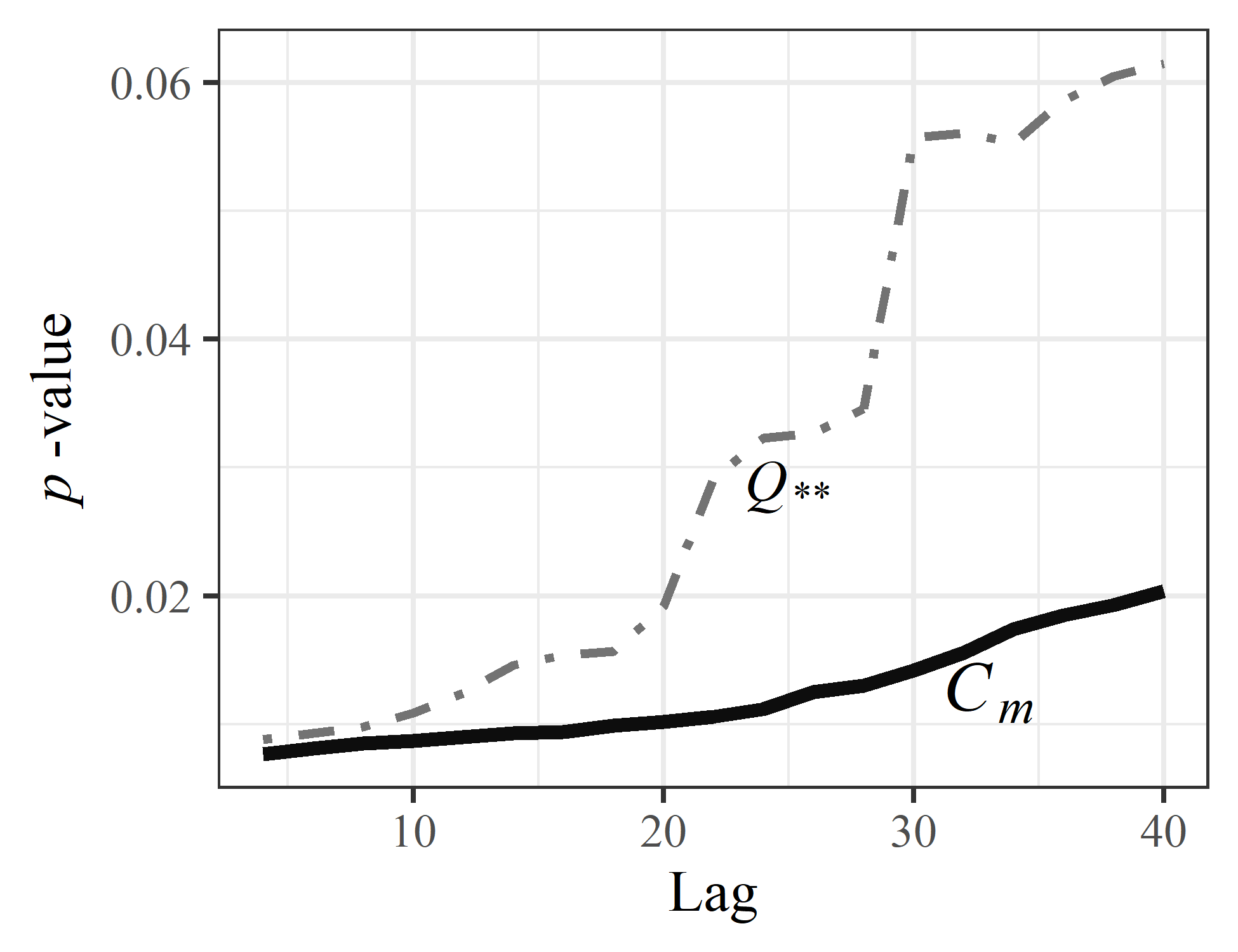}
	\end{center}
	\caption{{\color{blue}$p$-values of $C_m$ (solid black line) and $Q_{**}$ (dotted dashed gray line) at lags $m=4,\ldots, 40$ demonstrating the relative stability of the asymmetrically weighted statistic $C_m$ compared to $Q_{**}$ when an \AR(3) model is fit to the daily returns for the WTI Cude Oil Price series.}}
	\label{fig:wtiPvaluePlot}
\end{figure}
There, we see that the reported $p$-values of the asymmetrically weighted $C_m$ are relatively stable across all lags studied while the $Q_{**}$ has a noticeable increasing behavior demonstrating it is more sensitive to the chosen lag.

\subsection{Air Quality Measurements}
}
The Nitric Oxide measures (micrograms per cubic meter) at Marylebone Road and North Kensington air quality stations in London were collected from the Department for Environment, Food and Rural Affairs in the United Kingdom.
These two series are a subset of a multivariate time series explored in \cite{CirkovicFisher2021} and are available in their accompanying R package \texttt{autocovarianceTesting}.

\begin{figure}[htp!]
	\begin{center}
		\includegraphics[width=0.85\textwidth]{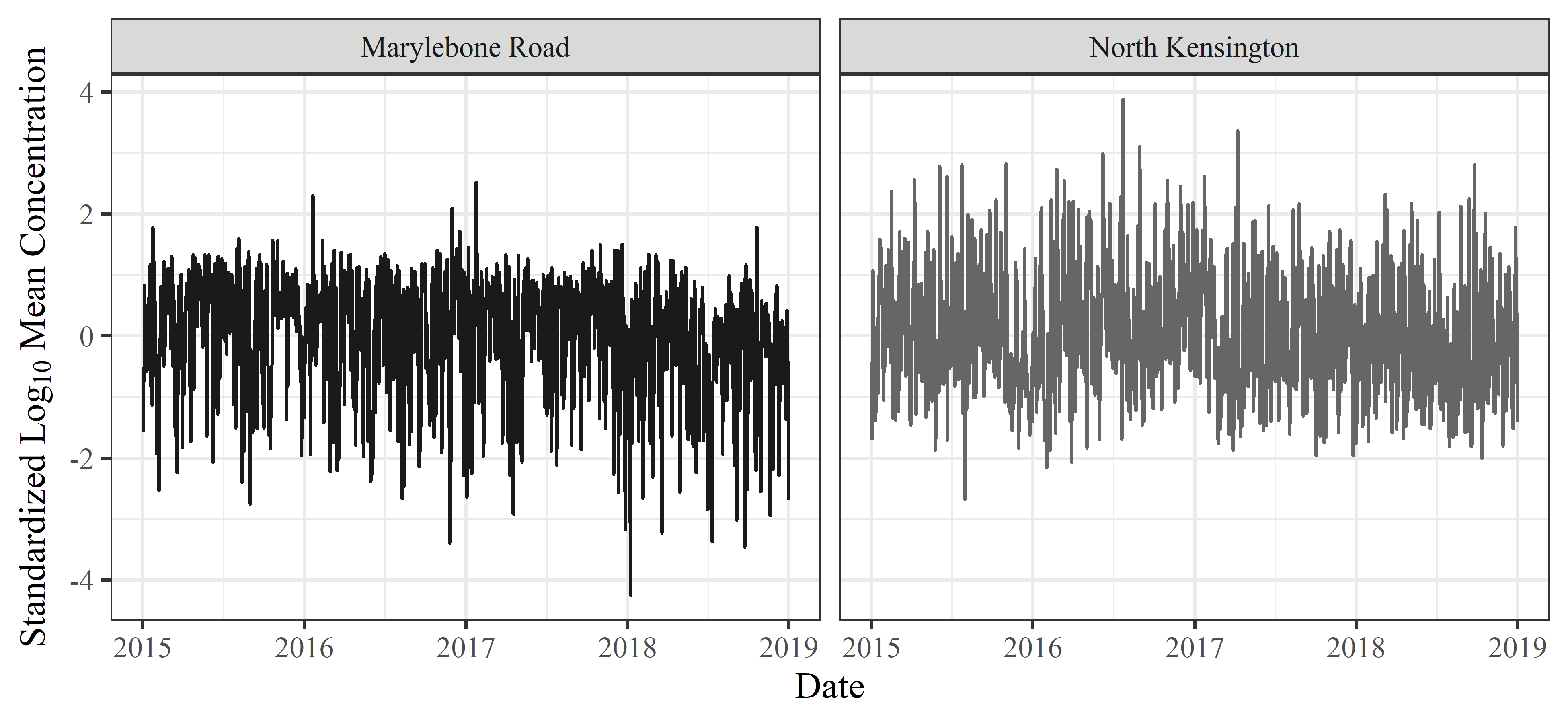}
	\end{center}
	\caption{Standardized logarithm (base 10) daily mean concentrations of Nitric Oxide at the Marylebone Road and North Kensington air quality stations in London, UK from 2015 to 2019.}
	\label{fig:londonAQ}
\end{figure}

Following the procedure in \cite{CirkovicFisher2021} the two $n=1461$ length marginal series are transformed by a logarithm and then standardized by both month and by weekday/weekend means and standard deviations to achieve stationary. 
The resulting series can be {\color{blue} seen} in Figure \ref{fig:londonAQ}.
The automatic lag selecting test procedure from \cite{CirkovicFisher2021} suggests the two series have equivalent autocovariance structures ($H_0:$\ two series share a common autocovariance structure, $p$-value of 0.504).
Thus, we may expect the two series to have similar linear dynamics and have a similar \ARMA{} fit.
In fact, using AIC to select the order, an \AR(3) is suggested for both series as an appropriate model.
The fitted \AR{} model parameters are {\color{blue}quite similar:}
\begin{align*}
	\textrm{Marylebone Road}:&~~ \phi_1=0.4542, ~\phi_2=0.0237, ~\phi_3=0.0435, ~~\sigma^2=0.7591 \\
	\textrm{North Kensington}:&~~ \phi_1=0.4252, ~\phi_2=0.0032, ~\phi_3=0.0503, ~~\sigma^2=0.7953
\end{align*}

\begin{figure}[htp!]
	\begin{center}
		\includegraphics[width=0.85\textwidth]{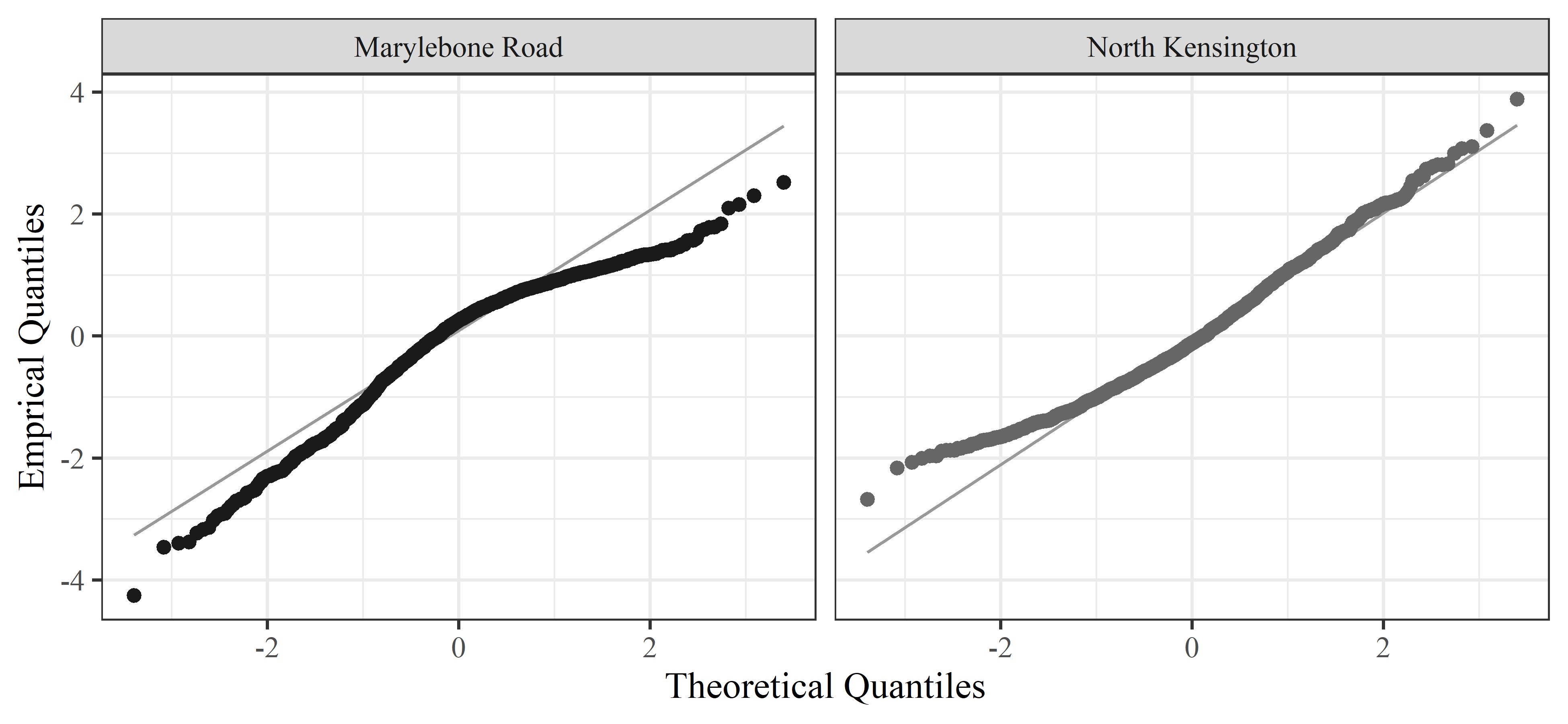}
	\end{center}
	\caption{Normal QQ-Plots of the standardized logarithm (base 10) daily mean concentrations of Nitric Oxide at the Marylebone Road and North Kensington Air Quality stations.}
	\label{fig:aqQQplot}
\end{figure}

As a follow-up, consider testing the adequacy of the fitted \AR(3) models for the two series by applying the proposed statistic and those studied from the literature.
The two series exhibit some {\color{blue} minor} deviations from normality (see Figure \ref{fig:aqQQplot}) so we utilize {\color{blue} both} the RWB algorithm to calculate the associated $p$-values of the test statistics as well as reporting the $p$-values based on the asymptotic distribution.
Table \ref{tab:aqPvals} reports the test statistic values and associated $p$-values for the six studied test statistics at two maximum lags.

\begin{table}
\tbl{The statistics and associated $p$-values, based on both the asymptotic distribution and RWB algorithm, of the portmanteau tests when an \AR(3) is fit to the two Nitric Oxide air quality datasets.}
{\begin{tabular}{llrrrccrrr}
					&& \multicolumn{3}{c}{Lag $m=5$} & & & \multicolumn{3}{c}{Lag $m=10$} \\
					\noalign{\smallskip}
					\cline{3-5}\cline{8-10}
					\noalign{\smallskip}
					& &Stat.\ & Asymp.\ & RWB & & & Stat.\ & Asymp.\ & RWB \\
					& &      & $p$-value & $p$-value & & & & $p$-value & $p$-value \\
					\hline\noalign{\smallskip}
					&&\multicolumn{8}{c}{\underline{Marylebone Road}}\\ 
					\noalign{\smallskip}
					$C_m$ & & 108.08 & $10^{-28}$ & $10^{-5}$ & & & 119.38 & $10^{-24}$ & $10^{-5}$ \\
					$Q_{**}$ & & 129.53 & $10^{-19}$ & $10^{-5}$ & & & 145.89 & $10^{-15}$ & $10^{-5}$ \\
					$Q_{11}$ & & 1.53 & 0.466 & 0.513 & & & 4.55 & 0.714 & 0.760 \\
					$Q_{22}$ & & 56.98 & $10^{-11}$ & $10^{-4}$ & & & 63.81 & $10^{-10}$ & $10^{-4}$ \\
					$Q_{12}$ & & 7.96 & 0.158 & 0.237 & & & 12.31 & 0.265 & 0.360 \\
					$Q_{21}$ & & 63.06 & $10^{-12}$ & $10^{-5}$ & & & 65.21 & $10^{-10}$ & $10^{-5}$ \\
					\hline\noalign{\smallskip}
					&&\multicolumn{8}{c}{\underline{North Kensington}}\\ 
					\noalign{\smallskip}
					$C_m$ & & 50.51 & $10^{-11}$ & $10^{-4}$ & & & 63.60 & $10^{-9}$ & $10^{-4}$ \\
					$Q_{**}$ & & 60.32 & $10^{-7}$ & 0.002 & & & 96.70 & $10^{-7}$ & $10^{-4}$ \\
					$Q_{11}$ & & 0.05 & 0.977 & 0.988 & & & 3.95 & 0.785 & 0.775 \\
					$Q_{22}$ & & 20.03 & 0.001 & 0.025 & & & 27.25 & 0.002 & 0.026 \\
					$Q_{12}$ & & 15.12 & 0.010 & 0.022 & & & 20.03 & 0.029 & 0.043 \\
					$Q_{21}$ & & 25.13 & $10^{-4}$ & 0.002 & & & 45.46 & $10^{-6}$ & $10^{-5}$ \\
					\hline
			\end{tabular}}
\label{tab:aqPvals}
\end{table}

\begin{figure}[htp!]
	\begin{center}
		\includegraphics[width=0.85\textwidth]{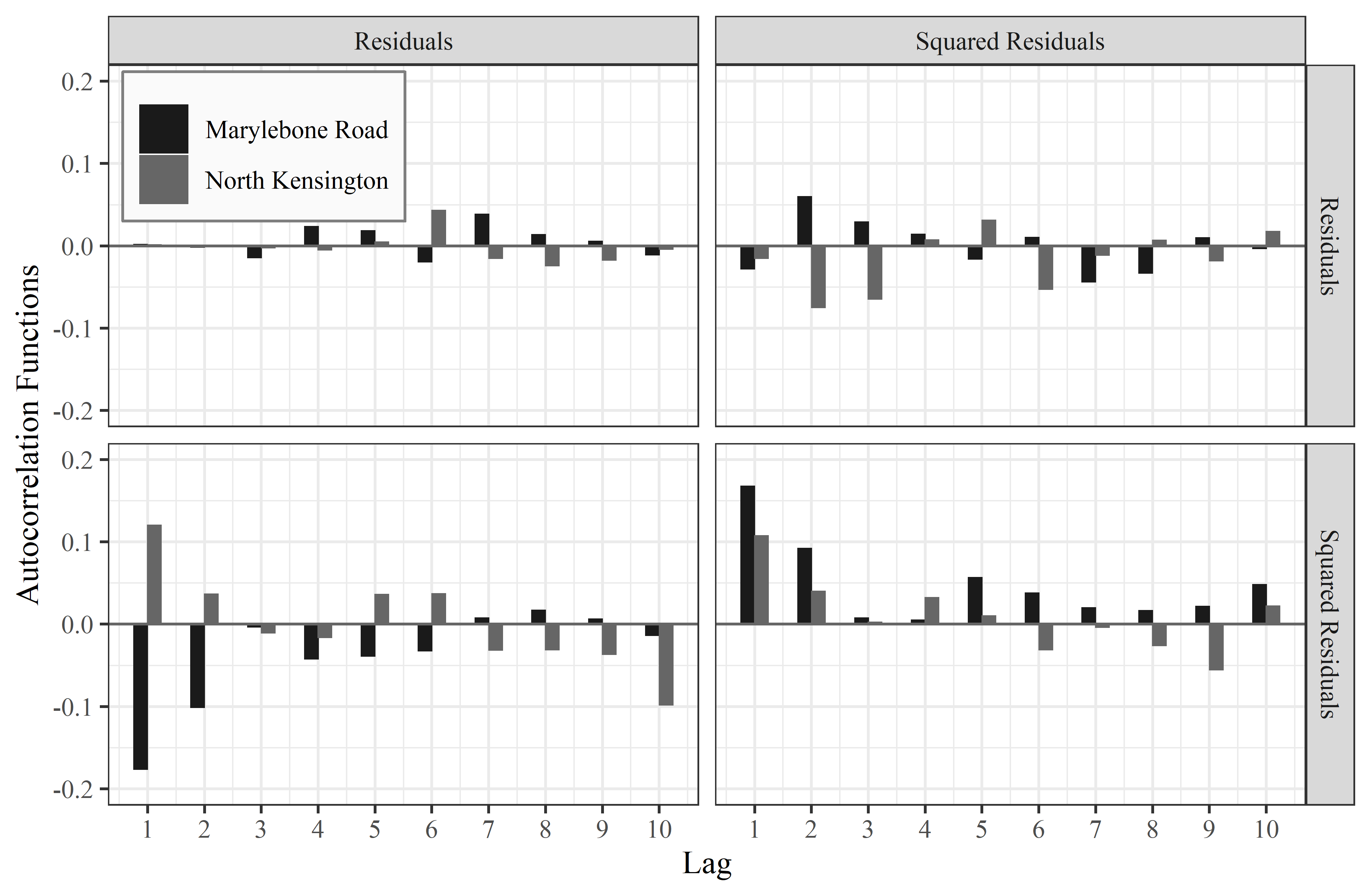}
	\end{center}
	\caption{Correlograms of the residuals and their squares from the two \AR(3) models for the standardized Marylebone Road and North Kensington Nitric Oxide series.}
	\label{fig:aqResACFs}
\end{figure}

We see that the two \AR(3) models adequately model the linear relationship based on the Ljung-Box statistic, $Q_{11}$, but there is overwhelming evidence for the presence of nonlinear effects. 
In particular, the proposed statistic $C_m$ provides unquestionable evidence for the inadequacy of the \AR(3) model as does the statistic of McLeod-Li, $Q_{22}$ and Pasaradakis-V\'{a}vra, $Q_{21}$.
To gain further insight, consider the correlogram plots of the two residual series and squared residuals in Figure \ref{fig:aqResACFs}.
The upper-left correlogram is a plot of the autocorrelation of the two residual series.
In agreement with the Ljung-Box statistic, there is no meaningful correlation present.
However, we see a decaying autocorrelation feature in the two bottom panels (corresponding to components of the $Q_{21}$ and $Q_{22}$ statistics). 
{\color{blue} Contextually, Figure \ref{fig:aqResACFs} appears to show that although the linear process of Marylebone Road and North Kensington Nitric Oxide series may be equivalent, the nonlinear processes differ (consider the bottom left panel of Figure \ref{fig:aqResACFs}).}

{\color{blue}
\subsection{Discussion on Applications}

We remind the reader the proposed statistic, $C_m$, is asymptotically equivalent to a convolution of $Q_{11}$, $Q_{22}$, $Q_{12}$ and $Q_{21}$, and that the statistic $Q_{**}$ is essentially the summation of the four.
In the WTI series, the statistics $Q_{12}$ and $Q_{22}$ consistently suggest nonlinearity, however with the air quality series, $Q_{21}$ and $Q_{22}$ suggest nonlinearity in the residuals while $Q_{12}$ provides contradictory results.
In both studies, the omnibus tests $C_m$ and $Q_{**}$ provide evidence for the presence of a nonlinear temporal structure.
Since foretelling which statistic is preferred for a given dataset would require oracle type abilities, these examples demonstrate the utility of using a omnibus statistic.

Further, the construction of $C_m$ follows that of \cite{PR2002} and has similarities to that of \cite{FisherGallagher2012}.
Those test essentially weight the correlation at lower lags with more emphasis than those at higher lags, while the test $Q_{**}$ considers all lags equally.
Likewise, the \cite{ZM2019} tests, $Q_{12}$ and $Q_{21}$, and \cite{McLeodLi1983} test, $Q_{22}$, also equally weigh each lag.
With the WTI data, Figure \ref{fig:wtiPvaluePlot} demonstrates the consistency of $C_m$ regardless of lag compared to $Q_{**}$.
For the air quailty data analysis in Figure \ref{fig:aqResACFs}, it appears the strongest evidence of nonlinearity is at lags 1 and 2, and we see that that $C_m$ offers more evidence than either $Q_{21}$, $Q_{22}$ and $Q_{**}$, particularly at the larger lag (Table \ref{tab:aqPvals}).

}

%%%%%%%%%%%%%%%%%%%%%%%%%%%%%%%%%
%%%%%%%%%%%%%%%%%%%%%%%%%%%%%%%%%
%% Discussion
%%%%%%%%%%%%%%%%%%%%%%%%%%%%%%%%%
%%%%%%%%%%%%%%%%%%%%%%%%%%%%%%%%%

%%%%%=============================================================================%%%%%
\section{\label{sec:discussion}Discussion}
%%%%%=============================================================================%%%%%
The proposed test statistic has several interesting properties. 
It can be seen as a combination of four weighted tests.
The first test is based on the partial autocorrelation of the residuals that can be used to test for linearity in time series models.
The second is based on the partial autocorrelation of the squared residuals that can be used to test for nonlinearity.  
The third and the fourth tests are based on the cross-correlation between the residuals and their squares at negative and positive lags, respectively.
Each term in the test is scaled by $(m-i+1)/m$, which allows the lower-order autocorrelations and cross-correlations to receive more emphasis than the larger lag terms.

In contrast to some other portmanteau tests, the proposed test responds well to nonlinear models that do not have \ARCH-type structures. 
In particular, the proposed test responds very well to time series where the residuals and their squares are cross correlated.  
Simulation results show the power of the proposed test is comparable, if not more powerful than, other nonlinear tests studied by \cite{McLeodLi1983} and \cite{ZM2019}.
% This test shows less sensitivity to Gaussian assumptions compared to the other portmanteau tests.

Several possible extensions to this article can be pursued. 
One is to approximate the distribution, and $p$-value of the proposed statistic, based on a different bootstrapping, or Monte Carlo, method than the RWB method used here.
For instance, Monte Carlo methods are suggested by \cite{LinMcLeod2006} and \cite{MahdiMcLeod2012} to compute $p$-values of a portmanteau test statistic based on the determinant autocorrelation matrix, but these methods are computationally expensive as they require repeated calculated of the determinant of a $2(m+1)\times 2(m+1)$ matrix required on the order of $O(4(m+1)^2)$ operations.

Another extension to this article could be done by deriving a new test based on extending the block matrix given by (\ref{eq:GCMat}) to the other generalized-correlation terms. 
Lastly, generalizing the result for the use of multivariate time series, similar to \cite{MahdiMcLeod2012} and \cite{RobbinsFisher2015} seems like a fairly straightforward calculation but may require very large samples to be tenable.

\bibliographystyle{tfs}
\bibliography{powerfulPortmeanteau}

\end{document}